\documentclass{m2an}
\usepackage{amsmath,amsfonts,amssymb,latexsym,epic,eepic}
\usepackage[dvips]{graphics,epsfig,color}

\newcommand{\dx}{\, {\rm d}x}
\newcommand{\dedge}{\, {\rm d} \gamma}
\newcommand{\dt}{\delta t}
\newcommand{\eqdef}{\stackrel{\mathrm{def}}{=}}
\newcommand{\vs}{\mathrm{v}_{\edge,K}}
\newcommand{\vsp}{\mathrm{v}_{\edge,K}^+}
\newcommand{\vsm}{\mathrm{v}_{\edge,K}^-}

\newcommand{\edge}{\sigma}
\newcommand{\edged}{\varepsilon}
\newcommand{\edges}{{\cal E}}
\newcommand{\edgesint}{{\cal E}_{{\rm int}}}
\newcommand{\edgesext}{{\cal E}_{{\rm ext}}}
\newcommand{\mesh}{{\cal M}}
\newcommand{\fluxK}{F_{\edge,K}}
\newcommand{\fluxL}{F_{\edge,L}}
\newcommand{\fluxd}{F_{\edged,\edge}}
\newcommand{\normLdiscd}[2]{\hspace{.2em}|\hspace{-.1em}|#2|\hspace{-.1em}|_{h,#1}^2\hspace{.2em}}
\newcommand{\snormundisc}[2]{|#2|_{h,#1}}
\newcommand{\snormundiscd}[2]{|#2|_{h,#1}^2}

\newcommand{\normLd}[1]{\hspace{.2em}|\hspace{-.1em}| #1 |\hspace{-.1em}|_{\xLtwo(\Omega)}\hspace{.2em}}

\newcommand{\normLdv}[1]{\hspace{.2em}|\hspace{-.1em}| #1 |\hspace{-.1em}|_{\xLtwo(\Omega)^d}\hspace{.2em}}

\newcommand{\normHb}[1]{\hspace{.2em}|\hspace{-.1em}| #1 |\hspace{-.1em}|_{1,b}\hspace{.2em}}
\newcommand{\normHbd}[1]{\hspace{.2em}|\hspace{-.1em}| #1 |\hspace{-.1em}|_{1,b}^2\hspace{.2em}}
\newcommand{\norms}[1]{\hspace{.2em}|\hspace{-.1em}| #1 |\hspace{-.1em}|_{\ast}\hspace{.2em}}
\newcommand{\normsd}[1]{\hspace{.2em}|\hspace{-.1em}| #1 |\hspace{-.1em}|_{\ast}^2\hspace{.2em}}
\newcommand{\normS}[1]{\hspace{.2em}|\hspace{-.1em}| #1 |\hspace{-.1em}|^{\ast}\hspace{.2em}}
\newcommand{\normSd}[1]{{\hspace{.2em}|\hspace{-.1em}| #1 |\hspace{-.1em}|^{\ast}}^2\hspace{.2em}}
\newcommand{\normz}[1]{\hspace{.2em}|\hspace{-.1em}| #1 |\hspace{-.1em}|\hspace{.2em}}
\newcommand{\matB}{{\rm B}}
\newcommand{\matL}{{\rm L}}
\newcommand{\matM}{{\rm M}}
\newcommand{\matR}{{\rm R}}
\newcommand{\matQ}{{\rm Q}}

\begin{document}

\title{An unconditionnally stable pressure correction scheme for compressible barotropic {N}avier-{S}tokes equations}

\author{T. Gallou\"et}
\address{Universit\'e de Provence, France (gallouet@cmi.univ-mrs.fr)}

\author{L. Gastaldo}
\address{Institut de Radioprotection et de S\^{u}ret\'{e} Nucl\'{e}aire (IRSN) (laura.gastaldo@irsn.fr)}

\author{R. Herbin}
\address{Universit\'e de Provence, France (herbin@cmi.univ-mrs.fr)}

\author{J.C. Latch\'e}
\address{Institut de Radioprotection et de S\^{u}ret\'{e} Nucl\'{e}aire (IRSN) (jean-claude.latche@irsn.fr)}

\begin{abstract}
We present in this paper a pressure correction scheme for barotropic compressible Navier-Stokes equations, which enjoys an unconditional stability property, in the sense that the energy and maximum-principle-based a priori estimates of the continuous problem also hold for the discrete solution.
The stability proof is based on two independent results for general finite volume discretizations, both interesting for their own sake: the $L^2$-stability of the discrete advection operator provided it is consistent, in some sense, with the mass balance and the estimate of the pressure work by means of the time derivative of the elastic potential.
The proposed scheme is built in order to match these theoretical results, and combines a fractional-step time discretization of pressure-correction type to a space discretization associating low order non-conforming mixed finite elements and finite volumes.
Numerical tests with an exact smooth solution show the convergence of the scheme.
\end{abstract}


%
%
%
%
%
%
%
\maketitle


\section{Introduction}

The problem addressed in this paper is the system of the so-called barotropic compressible Navier-Stokes equations, which reads:
\begin{equation}
\left| \begin{array}{l} \displaystyle
\frac{\partial\,\rho}{\partial t}+\nabla\cdot(\rho\, u)=0 
\\[2ex] \displaystyle
\frac{\partial}{\partial t}(\rho\, u)+\nabla\cdot(\rho\, u\otimes u)+\nabla p -\nabla\cdot\tau(u)=f 
\\[2ex] \displaystyle
\rho=\varrho(p)
\end{array} \right.
\label{nsb}\end{equation}
where $t$ stands for the time, $\rho$, $u$ and $p$ are the density, velocity and pressure in the flow, $f$ is a forcing term and $\tau(u)$ stands for the shear stress tensor.
The function $\varrho(\cdot)$ is the equation of state used for the modelling of the particular flow at hand, which may be the actual equation of state of the fluid or may result from assumptions concerning the flow; typically, laws as $\rho=p^{1/\gamma}$, where $\gamma$ is a coefficient specific to the fluid considered, are obtained by making the assumption that the flow is isentropic.
This system of equations is posed over $\Omega\times (0,T)$, where $\Omega$ is a domain of $\xR^d$, $d\leq 3$ supposed to be polygonal ($d=2$) or polyhedral ($d=3$), and the final time $T$ is finite.
It must be supplemented by boundary conditions and by an initial condition for $\rho$ and $u$.

\medskip
The development of pressure correction techniques for compressible Navier-Stokes equations may be traced back to the seminal work of Harlow and Amsden \cite{har-68-num, har-71-num} in the late sixties, who developped an iterative algorithm (the so-called ICE method) including an elliptic corrector step for the pressure.
Later on, pressure correction equations appeared in numerical schemes proposed by several researchers, essentially in the finite-volume framework, using either a collocated \cite{pat-87-bar,dem-93-col,kob-96-cha,pol-97-pre,iss-98-pre,mou-01-hig} or a staggered arrangement \cite{cas-84-pre, iss-85-sol, iss-86-com, kar-89-pre, bij-98-uni, col-99-pro, van-01-sta, wal-02-sem, wen-02-mac, van-03-con, vid-06-sup} of unknowns; in the first case, some corrective actions are to be foreseen to avoid the usual odd-even decoupling of the pressure in the low Mach number regime.
Some of these algorithms are essentially implicit, since the final stage of a time step involves the unknown at the end-of-step time level; the end-of-step solution is then obtained by SIMPLE-like iterative processes \cite{van-87-seg, kar-89-pre, dem-93-col, kob-96-cha, pol-97-pre, iss-98-pre, mou-01-hig}.
The other schemes \cite{iss-85-sol, iss-86-com, pat-87-bar, bij-98-uni, col-99-pro, wes-01-pri, van-01-sta, wen-02-mac, van-03-con, vid-06-sup} are predictor-corrector methods, where basically two steps are performed in sequence: first a semi-explicit decoupled prediction of the momentum or velocity (and possibly energy, for non-barotropic flows) and, second, a correction step where the end-of step pressure is evaluated and the momentum and velocity are corrected, as in projection methods for incompressible flows (see \cite{cho-68-num, tem-69-sur} for the original papers, \cite{mar-98-nav} for a comprehensive introduction and \cite{gue-06-ove} for a review of most variants).
The Characteristic-Based Split (CBS) scheme (see \cite{nit-06-cha} for a recent review or \cite{zie-95-gen} for the seminal paper), developped in the finite-element context, belongs to this latter class of methods.

\medskip
Our aim in this paper is to propose and study a simple non-iterative pressure correction scheme for the solution of \eqref{nsb}.
In addition, this method is designed so as to be stable in the low Mach number limit, since our final goal is to apply it to simulate by a drift-flux approach a class of bubbly flows encountered in nuclear safety studies, where pure liquid (incompressible) and pure gaseous (compressible) zones may coexist.
To this purpose, we use a low order mixed finite element approximation, which meets the two following requirements: to allow a natural discretization of the viscous terms and to provide a spatial discretization that is intrinsically stable (\ie \ without the adjunction of stabilization terms to circumvent the so-called {\it inf-sup} or BB condition) in the incompressible limit.

\medskip
In this work, a special attention is payed to stability issues.
To be more specific, let us recall the {\it a priori} estimates associated to problem \eqref{nsb} with a zero forcing term, {\it i.e.} estimates which should be satisfied by any possible regular solution \cite{pll-98-mat, fei-04-dyn, nov-04-int}:
\begin{equation}
\left| \quad \begin{array}{llll}
(i)
& \displaystyle
\rho(x,t) > 0,
&& \displaystyle
\forall x \in \Omega,\ \forall t \in (0,T)
\\[2ex]
(ii)
&  \displaystyle
\int_\Omega \rho(x,t)\ {\rm d} x = \int_\Omega \rho(x,0)\ {\rm d} x, 
&& \displaystyle
\forall t \in (0,T) 
\\[3ex]
(iii) \quad
& \displaystyle
\frac{1}{2}\ \frac{d}{dt}\int_{\Omega}\rho(x,t)\, u(x,t)^{2} \dx +\frac{d}{dt}\int_{\Omega} \rho(x,t) \,P(\rho(x,t)) \dx 
\\
& \hfill \displaystyle
+\int_{\Omega} \tau (u(x,t)):\nabla u(x,t) \dx = 0,
& \qquad &
\forall t \in (0,T)
\end{array} \right.
\label{apriori-e}\end{equation}
In the latter relation, $P$, the "elastic potential", is a function derived from the equation of state, which satisfies:
\begin{equation}
P'(z)=\frac{\wp(z)}{z^2}
\label{elasticpot_gene} \end{equation}
where $\wp (\cdot)$ is the inverse function of $\varrho(\cdot)$, \ie\ the function giving the pressure as a function of the density.
The usual choice is, provided that this expression makes sense:
\begin{equation}
P(z)=\int_0^z \frac{\wp(s)}{s^2}\ {\rm d}s
\label{elasticpot} \end{equation}
For these estimates to hold, the condition \eqref{apriori-e}-$(i)$ must be satisfied by the initial condition; note that a non-zero forcing term $f$ in the momentum balance would make an additional term appear at the right hand side of relation \eqref{apriori-e}-$(iii)$.
This latter estimate is obtained from the second relation of \eqref{nsb} (\ie\ the momentum balance) by taking the scalar product by $u$ and integrating over $\Omega$.
This computation then involves two main arguments which read:
\begin{equation}
\begin{array}{llll}
(i)
&
\mbox{Stability of the advection operator:} \quad
& \displaystyle
\int_\Omega \left[ \frac{\partial}{\partial t}(\rho\, u)+\nabla\cdot(\rho\, u\otimes u) \right] \cdot u \dx=
\frac{1}{2}\ \frac{d}{dt}\int_{\Omega}\rho\, u^2 \dx
\\[3ex]
(ii)
&
\mbox{Stability due to the pressure work:} \quad
& \displaystyle
\int_\Omega -p \, \nabla \cdot u \dx=\frac{d}{dt}\int_{\Omega} \rho(x,t) \,P(\rho(x,t)) \dx
\end{array}
\label{stab-arg}\end{equation}
Note that the derivation of both relations make use of the fact that the mass balance equation holds in a crucial way.

\medskip
This paper is organized as follows.

As a first step, we derive a bound similar to \eqref{stab-arg}-$(ii)$ for a given class of spatial discretizations; the latter are introduced in section \ref{subsec:disc} and the desired stability estimate (theorem \ref{VF2}) is stated and proven in section \ref{subsec:pwork}.
We then show that this result allows to prove the existence of a solution for a fairly general class of discrete compressible flow problems.
Section \ref{sec:darcy} gathers this whole study, and constitutes the first part of this paper.

In a second part (section \ref{sec:scheme}), we turn to the derivation of a pressure correction scheme the solution of which satisfies the whole set of {\it a priori} estimates \eqref{apriori-e}.
To this purpose, besides theorem \ref{VF2}, we need as a second key ingredient a discrete version of the bound \eqref{stab-arg}-$(i)$ relative to the stability of the advection operator, which is stated and proven in section \ref{subsec:conv} (theorem \ref{VF1}).
We then derive a fully discrete algorithm which is designed to meet the assumptions of these theoretical results, and establish its stability.
Some numerical experiments show that this scheme seems in addition to present, when the solution is smooth, the convergence properties which can be expected, namely first order in time convergence for all the variables and second order in space in $\xLtwo$ and discrete $\xLtwo$ norm for the velocity and the pressure, respectively.


\section{Analysis of a class of discrete problems} \label{sec:darcy}

The class of problems addressed in this section can be seen as the class of discrete systems obtained by space discretization by low-order non-conforming finite elements of continuous problems of the following form:
\begin{equation}
\left| \begin{array}{ll} \displaystyle
A \, u + \nabla p =  f
& \mbox{ in } \Omega
\\[2ex] \displaystyle
\frac{\varrho(p)-\rho^\ast}{\dt}+ \nabla \cdot \left( \varrho(p)\, u \right)=0
& \mbox{ in } \Omega
\\[2ex] \displaystyle
u=0 & \mbox{ on } \partial \Omega
\end{array} \right .
\label{pbdisc-cont}\end{equation}
where the forcing term $f$ and the density field $\rho^\ast$ are known quantities, and $A$ stands for an abstract elliptic operator.
The unknown of the problem are the velocity $u$ and the pressure $p$; the function $\varrho(\cdot)$ stands for the equation of state.
The domain $\Omega$ is a polygonal ($d=2$) or polyhedral ($d=3$) open, bounded and connected subset of $\xR^d$.
Of course, at the continuous level, this statement of the problem should be completed by a precise definition of the functional spaces in which the velocity and the pressure are searched for, together with regularity assumptions on the data.
This is out of context in this section, as system \eqref{pbdisc-cont} is given only to fix ideas and we restrict ourselves here to prove some mathematical properties of the discrete problem, namely to establish some {\it a priori} estimates for its solution and to prove that this nonlinear problem admits some solution for fairly general equations of state.

\medskip
This section is organized as follows.
We begin by describing the considered discretization and precisely stating the discrete problem at hand.
Then we prove, for the particular discretization at hand, a fundamental result which is a discrete analogue of the elastic potential identity.
The next section is devoted to the proof of the existence of a solution, and we finally conclude by giving some practical examples of application of the abstract theory developped here.


\subsection{The discrete problem}\label{subsec:disc}

Let $\mesh$ be a decomposition of the domain $\Omega$ either into convex quadrilaterals ($d=2$) or hexahedrons ($d=3$) or in simplices. 
By $\edges$ and $\edges(K)$ we denote the set of all $(d-1)$-edges $\edge$ of the mesh and of the element $K \in \mesh$ respectively.
The set of edges included in the boundary of $\Omega$ is denoted by $\edgesext$ and the set of internal ones (\ie\ $\edges \setminus \edgesext$) is denoted by $\edgesint$.
The decomposition $\mesh$ is supposed to be regular in the usual sense of the finite element literature (e.g. \cite{cia-91-bas}), and, in particular, $\mesh$ satisfies the following properties: $ \bar\Omega=\bigcup_{K\in \mesh} \bar K$; if $K,\,L \in \mesh,$ then $\bar K \cap \bar L=\emptyset$ or $\bar K\cap \bar L$ is a common face of $K$ and $L$, which is denoted by $K|L$.
For each internal edge of the mesh $\edge=K|L$, $n_{KL}$ stands for the normal vector of $\edge$, oriented from $K$ to $L$.
By $|K|$ and $|\edge|$ we denote the measure, respectively, of $K$ and of the edge $\edge$.

\medskip
The spatial discretization relies either on the so-called "rotated bilinear element"/$P_0$ introduced by Rannacher and Turek \cite{ran-92-sim} for quadrilateral of hexahedric meshes, or on the Crouzeix-Raviart element (see \cite{cro-73-con} for the seminal paper and, for instance, \cite[p. 83--85]{ern-05-aid} for a synthetic presentation) for simplicial meshes.
The reference element $\widehat K$ for the rotated bilinear element is the unit $d$-cube (with edges parallel to the coordinate axes); the discrete functional space on $\widehat K$ is $\tilde{Q}_{1}(\widehat K)^d$, where $\tilde{Q}_{1}(\widehat K)$ is defined as follows:
\[
\tilde{Q}_{1}(\widehat K)= {\rm span}\left\{1,\,(x_{i})_{i=1,\ldots,d},\,(x_{i}^{2}-x_{i+1}^{2})_{i=1,\ldots,d-1}\right\}
\]
The reference element for the Crouzeix-Raviart is the unit $d$-simplex and the discrete functional space is the space $P_1$ of affine polynomials.
For both velocity elements used here, the degrees of freedom are determined by the following set of nodal functionals:
\begin{equation}
\displaystyle \left\{F_{\edge,i},\ \edge \in \edges(K),\, i=1,\ldots,d\right\}, 
\qquad F_{\edge,i}(v)=|\edge|^{-1}\int_{\edge} v_{i} \dedge
\label{vdof}\end{equation}
The mapping from the reference element to the actual one is, for the Rannacher-Turek element, the standard $Q_1$ mapping and, for the Crouzeix-Raviart element, the standard affine mapping.
Finally, in both cases, the continuity of the average value of discrete velocities (\ie, for a discrete velocity field $v$, $F_{\sigma,i}(v)$) across each face of the mesh is required, thus the discrete space $W_{h}$ is defined as follows:
\[
\begin{array}{ll} \displaystyle
W_h =
& \displaystyle
\lbrace
\ v_h\in L^{2}(\Omega)\,:\, v_h|_K \in\tilde{Q}_{1}(K)^d,\,\forall K\in \mesh;\ F_{\edge,i}(v_h) \mbox{ continuous across each edge } \sigma \in {\edgesint},\ 1\leq i \leq d \,;
\\[1ex] & \displaystyle
\ \ F_{\edge,i}(v_h)=0,\ \forall \edge \in \edgesext,\ 1\leq i \leq d\ \rbrace
\end{array}
\]
For both Rannacher-Turek and Crouzeix-Raviart discretizations, the pressure is approximated by the space $L_{h}$ of piecewise constant functions:
\[
L_h=\left\{q_h\in L^{2}(\Omega)\,:\, q_h|_K=\mbox{ constant},\,\forall K\in \mesh\right\}
\]
Since only the continuity of the integral over each edge of the mesh is imposed, the velocities are discontinuous through each edge; the discretization is thus nonconforming in $H^1(\Omega)^d$.
These pairs of approximation spaces for the velocity and the pressure are \textit{inf-sup} stable, in the usual sense for "piecewise $\xHone$" discrete velocities, \ie\ there exists $c_{\rm i}>0$ independent of the mesh such that:
\[
\forall p \in L_h, \qquad \sup_{v \in W_h} \frac{\displaystyle \int_{\Omega,h} p \, \nabla \cdot v \dx}{\normHb{v}} \geq c_{\rm i} \normLd{p-\bar p}
\]
where $\bar p$ is the mean value of $p$ over $\Omega$, the symbol $\displaystyle \int_{\Omega,h}$ stands for $\displaystyle \sum_{K\in\mesh} \int_K$ and $\normHb{\cdot}$ stands for the broken Sobolev $\xHone$ semi-norm:
\[
\normHbd{v}=\sum_{K\in \mesh} \int_K |\nabla v |^2 \dx=\int_{\Omega,h}| \nabla v |^2 \dx
\]

\bigskip
From the definition \eqref{vdof}, each velocity degree of freedom can be univoquely associated to an element edge.
We take benefit of this correspondence by using hereafter, somewhat improperly, the expression "velocity on the edge $\edge$" to name the velocity vector defined by the  degrees of freedom of the velocity components associated to $\edge$.
In addition, the velocity degrees of freedom are indexed by the number of the component and the associated edge, thus the set of velocity degrees of freedom reads:
\[
\lbrace v_{\edge, i},\ \edge \in \edgesint,\ 1 \leq i \leq d \rbrace
\]
We define $v_\edge=\sum_{i=1}^d v_{\edge, i} e_i $ where $e_i$ is the $i^{th}$ vector of the canonical basis of $\xR^d$.
We denote by $\varphi_\edge^{(i)}$ the vector shape function associated to $v_{\edge, i}$, which, by the definition of the considered finite elements, reads:
\[
\varphi_\edge^{(i)}=\varphi_\edge \, e_i
\]
where $\varphi_\edge$ is a scalar function.
Similarly, each degree of freedom for the pressure is associated to a mesh $K$, and the set of pressure degrees of freedom is denoted by
$\lbrace p_K,\ K \in \mesh \rbrace$.

\bigskip
For any $K \in \mesh$, let $\rho^\ast_K$ be a quantity approximating a known density $\rho^\ast$ on $K$.
The family of real numbers $(\rho^\ast_K)_{K\in\mesh}$ is supposed to be positive.
The discrete problem considered in this section reads:
\begin{equation}
\left| \begin{array}{ll} \displaystyle
a(u,\varphi_\edge^{(i)}) - \int_{\Omega,h} p\ \nabla \cdot \varphi_\edge^{(i)} = \int_\Omega f \cdot \varphi_\edge^{(i)} \dx
\qquad & \displaystyle
\forall \edge\in\edgesint \ (\edge=K|L),\ 1 \leq i \leq d
\\[2ex] \displaystyle
|K|\ \frac{\varrho(p_K)-\rho^\ast_K}{\dt} + \sum_{\edge=K|L} \vsp\, \varrho(p_K) - \vsm\ \varrho(p_L)=0
\qquad & \displaystyle
\forall K \in \mesh
\end{array}\right .
\label{pbdisc-disc}\end{equation}
where $\vsp$ and $\vsm$ stands respectively for $\vsp = \max (\vs,\ 0)$ and $\vsm = \max ( -\vs,\ 0)$ with $\vs=|\edge|\, u_\edge \cdot n_{KL}=\vsp- \vsm$.
The first equation is the standard finite element discretization of the first relation of \eqref{pbdisc-cont}, provided the following identity holds:
\[
\forall v\in W_h,\ \forall w \in W_h,\qquad a(v,w)=\int_\Omega Av \cdot w \dx
\]
As the pressure is piecewise constant, the finite element discretization of the second relation of \eqref{pbdisc-cont}, \ie\ the mass balance, is similar to a finite volume formulation, in which we introduce the standard first-order upwinding.
The bilinear form $a(\cdot,\cdot)$ is supposed to be elliptic on $W_h$, \ie\ to be such that the following property holds:
\[
\exists c_{\rm a} >0 \mbox{ such that, } \forall v \in W_h,\qquad a(u,u) \geq c_{\rm a} \normsd{u}
\]
where $\norms{\cdot}$ is a norm over $W_h$.
We denote by $\normS{\cdot}$ its dual norm with respect to the $\xLtwo(\Omega)^d$ inner product, defined by:
\[
\forall v \in W_h, \qquad \normS{v}=\sup_{w \in W_h} \frac{\displaystyle \int_\Omega v \cdot w \dx}{\norms{w}}
\]


\subsection{On the pressure control induced by the pressure forces work}\label{subsec:pwork}

The aim of this subsection is to prove that the discretization at hand satisfies a stability bound which can be seen as the discrete analogue of equation \eqref{stab-arg}-$(ii)$, which we recall here:
\[
\int_\Omega p \, \nabla \cdot u \dx=\frac{d}{dt}\int_{\Omega} \rho \,P(\rho) \dx,
\qquad \mbox{where } P(\cdot) \mbox{ is such that } P'(z)=\frac{\wp(z)}{z^2}
\]
The formal computation which allows to derive this estimate in the continuous setting is the following.
The starting point is the mass balance, which is multiplied by $[\rho P(\rho)]'$:
\[
[\rho P(\rho)]' \left[ \frac{\partial \rho}{\partial t} + \nabla \cdot (\rho u) \right]=0
\]
This relation yields:
\begin{equation}
\frac{\partial [\rho P(\rho)]}{\partial t} + [\rho P(\rho)]' \left[ u \cdot \nabla \rho + \rho \nabla \cdot u \right]=0
\label{elP-step1}\end{equation}
And thus:
\[
\frac{\partial [\rho P(\rho)]}{\partial t} + u \cdot \nabla [\rho P(\rho)] + [\rho P(\rho)]' \rho \nabla \cdot u =0
\]
Developping the derivative, we get:
\begin{equation}
\frac{\partial [\rho P(\rho)]}{\partial t} + 
\underbrace{u \cdot \nabla [\rho P(\rho)] + \rho P(\rho) \nabla \cdot u}_{\displaystyle \nabla \cdot (\rho P(\rho)\, u)} +
\underbrace{\rho^2 [P(\rho)]' \nabla \cdot u}_{\displaystyle p \nabla \cdot u} =0
\label{elP-step2}\end{equation}
and the result follows by integration in space.
We are going here to reproduce this computation at the discrete level.
%
%
\begin{thrm}[Stability due to the pressure work]\label{VF2}
Let the elastic potential $P$ defined by \eqref{elasticpot_gene} be such that the function $z \mapsto z\, P(z)$ is once continuously differentiable and strictly convex.
Let $(\rho_K)_{K\in\mesh}$ and $(p_K)_{K\in\mesh}$ satisfy the second relation of \eqref{pbdisc-disc}.
Then the following estimate holds:
\begin{equation}
- \int_\Omega p \, \nabla \cdot u \dx=\sum_{K\in\mesh} - p_K  \sum_{\edge=K|L} \vs \geq
\frac{1}{\dt}\ \sum_{K\in\mesh}  |K| \left[ \rho_K \, P(\rho_K) - \rho_K^\ast \, P(\rho_K^\ast)\right]
\label{pot-el}\end{equation}
\end{thrm}
\begin{proof}
Let us write the divergence term in the discrete mass balance over $K$ (\ie\ the second relation of \eqref{pbdisc-disc}) under the following form:
\[
\sum_{\edge=K|L} \rho_\edge \, \vs
\]
where $\rho_\edge$ is either $\rho_K=\varrho(p_K)$ if $\vs \geq 0$ or $\rho_L=\varrho(P_L)$ if $\vs \leq 0$.
Multiplying this term by the derivative with respect to $\rho$ of $\rho\, P(\rho)$ computed at $\rho_K$, denoted by $[\rho\, P(\rho)]'_{\rho_K}$, we obtain:
\[
T_{{\rm div},K}= [\rho\, P(\rho)]'_{\rho_K} \sum_{\edge=K|L} \rho_\edge \, \vs
= [\rho\, P(\rho)]'_{\rho_K} \left[ \sum_{\edge=K|L} \ (\rho_\edge-\rho_K) \, \vs+
\rho_K \sum_{\edge=K|L} \vs \right]
\]
This relation is a discrete equivalent to equation \eqref{elP-step1}: up to the multiplication by $1/|K|$, the first summation in the right hand side is the analogue of $u \cdot \nabla \rho$ and the second one to $\rho\,\nabla \cdot u$.
Developping the derivative, we obtain the equivalent of relation \eqref{elP-step2}:
\begin{equation}
T_{{\rm div},K}=
[\rho\, P(\rho)]'_{\rho_K} \sum_{\edge=K|L} (\rho_\edge-\rho_K) \, \vs
+ \rho_K\, P(\rho_K) \sum_{\edge=K|L} \vs
+ \rho_K^2\, P'(\rho_K) \sum_{\edge=K|L} \vs
\label{elPd1}\end{equation}
By definition \eqref{elasticpot_gene} of $P$, the last term reads $p_K \sum_{\edge=K|L} \vs$.
The process will be completed if we put the first two terms in divergence form.
To this end, let us sum up the $T_{{\rm div},K}$ and reorder the summation:
\begin{equation}
\sum_{K\in\mesh} T_{{\rm div},K}= \sum_{K\in\mesh} p_K \sum_{\edge=K|L} \vs
+ \sum_{\edge\in\edgesint} T_{{\rm div},\edge}
\label{elPd2}\end{equation}
where, if $\edge=K|L$:
\[
T_{{\rm div},\edge}= \vs\  \left[ \rho_K\, P(\rho_K) + [\rho\, P(\rho)]'_{\rho_K} (\rho_\edge-\rho_K)
- \rho_L\, P(\rho_L) - [\rho\, P(\rho)]'_{\rho_L} (\rho_\edge-\rho_L) \right]
\]
In this relation, there are two possible choices for the orientation of $\edge$, \ie\ $K|L$ or $L|K$, where $K$ and $L$ are two cells such that $\edge=\bar K\cap \bar L$; we choose  this orientation to have $\vs \geq 0$.
Let $\bar \rho_\edge$ be defined by:
\begin{equation}
\left| \begin{array}{ll}
\mbox{if }\rho_K \neq \rho_L \ :\qquad
& \displaystyle
\rho_K\, P(\rho_K) + [\rho\, P(\rho)]'_{\rho_K} (\bar \rho_\edge-\rho_K) = \rho_L\, P(\rho_L) + [\rho\, P(\rho)]'_{\rho_L} (\bar \rho_\edge-\rho_L)
\\[2ex]
\mbox{otherwise}:
&
\bar \rho_\edge=\rho_K=\rho_L
\end{array} \right.
\label{defrhob}\end{equation}
As the function $z \mapsto z\, P(z)$ is supposed to be once continuously differentiable and strictly convex, the technical lemma \ref{rhosigma} proven hereafter applies and $\bar \rho_\edge$ is uniquely defined and satisfies $\bar\rho_{\edge} \in [\min(\rho_{K},\rho_{L}), \, \max(\rho_{K},\rho_{L})]$.
By definition, the choice $\rho_\edge=\bar \rho_\edge$ is such that the term $T_{{\rm div},\edge}$ vanishes, which means that, by this definition, we would indeed have obtained that the first two terms of equation \eqref{elPd1} are a conservative approximation of the quantity $\nabla \cdot \rho P(\rho) u$ appearing in equation \eqref{elP-step2}, with the following expression for the flux:
\[
\begin{array}{l} \displaystyle
F_{\edge,K}= [\rho P(\rho)]_\edge \, \vs, \quad \mbox{with:}
\\[2ex] \displaystyle \hspace{10ex}
[\rho P(\rho)]_\edge= \rho_K\, P(\rho_K) + [\rho\, P(\rho)]'_{\rho_K} (\bar \rho_\edge-\rho_K) = \rho_L\, P(\rho_L) + [\rho\, P(\rho)]'_{\rho_L} (\bar \rho_\edge-\rho_L)
\end{array}
\]
Whatever the choice for $\rho_\edge$ may be, we have:
\[
T_{{\rm div},\edge}= \vs\ (\rho_\edge - \bar \rho_\edge)\ ([\rho\, P(\rho)]'_{\rho_K}-[\rho\, P(\rho)]'_{\rho_L} )
\]
With the orientation taken for $\edge$, an upwind choice yields:
\[
T_{{\rm div},\edge}= \vs\ (\rho_K - \bar \rho_\edge)\ ([\rho\, P(\rho)]'_{\rho_K}-[\rho\, P(\rho)]'_{\rho_L} )
\]
and, using the fact that $z \mapsto [\rho\, P(\rho)]'_{z}$ is an increasing fuction since $z \mapsto z\, P(z)$ is convex and that $\bar\rho_{\edge} \in [\min(\rho_{K},\rho_{L}), \, \max(\rho_{K},\rho_{L})]$, it is easily seen that $T_{{\rm div},\edge}$ is non-negative.

\medskip
Multiplying by $[\rho\, P(\rho)]'_{\rho_K}$ the mass balance over each cell $K$ and summing for $K\in\mesh$ thus yields, by equation \eqref{elPd2}:
\begin{equation}
- \sum_{K\in\mesh} p_K \sum_{\edge=K|L} \vs = 
R + \sum_{K\in\mesh} |K|\ [\rho\, P(\rho)]'_{\rho_K}\ \frac{\rho_K-\rho^\ast}{\dt}
\label{elPd3}\end{equation}
where $R$ is non-negative, and the result follows invoking once again the convexity of $z \mapsto z\, P(z)$.
\end{proof}
\begin{rmrk}[On a non-dissipative scheme]
The preceding proof shows that, for a scheme to conserve the energy (\ie\ to obtain a discrete equivalent of \eqref{apriori-e}-$(iii)$), besides other arguments, the choice of $\bar \rho_\edge$ given by \eqref{defrhob} for the density at the face of the control volume in the discretization of the flux in the mass balance seems to be mandatory; any other choice leads to an artificial dissipation in the work of the pressure forces.
Note however that, this discretization being essentially of centered type, the positivity of the density is not warranted in this case.
\end{rmrk}
%
%
\bigskip
In the course of the preceding proof, we used the following technical lemma.
\begin{lmm}\label{rhosigma} Let $g(\cdot)$ be a strictly convex and once continuously derivable real function.
Let $\edge$ be an internal edge of the mesh separating the cells $K$ and $L$.
Then the following relations:
\begin{equation}
\left| \begin{array}{ll}
\mbox{if }\rho_K \neq \rho_L \ :\qquad
& \displaystyle
g(\rho_{K}) + g'(\rho_{K})(\bar\rho_\edge-\rho_{K}) = g(\rho_{L}) + g'(\rho_{L})(\bar\rho_\edge-\rho_{L})
\\[2ex]
\mbox{otherwise}:
&
\bar \rho_\edge=\rho_K=\rho_L
\end{array} \right.
\label{defrhob2}\end{equation}
uniquely defines the real number $\bar\rho_\edge$.
In addition, we have $\bar\rho_{\edge} \in [\min(\rho_{K},\rho_{L}),\,\max(\rho_{K},\rho_{L})]$.
\end{lmm}

\begin{proof}
If $\rho_K=\rho_L$, there is nothing to prove.
In the contrary case, without loss of generality, let us choose $K$ and $L$ in such a way that $\rho_K < \rho_L$.
By reordering equation \eqref{defrhob2}, we get:
\[
g(\rho_K) + g'(\rho_K)\,(\rho_L-\rho_K) -g(\rho_L)= (\bar \rho_\edge-\rho_L)\,\left[ g'(\rho_L)-g'(\rho_K)\right]
\]
As, because $g(\cdot)$ is strictly convex, $g'(\rho_L)-g'(\rho_K)$ does not vanish, this equation proves that $\bar\rho_\edge$ is uniquely defined.
In addition, for the same reason, the left hand side of this relation is negative and $g'(\rho_L)-g'(\rho_K)$ is positive, thus we have $\bar \rho_\edge < \rho_L$.
Reordering in another way equation \eqref{defrhob2} yields:
\[
g(\rho_L) + g'(\rho_L)\,(\rho_K-\rho_L) -g(\rho_K)= (\bar \rho_\edge-\rho_K)\,\left[ g'(\rho_K)-g'(\rho_L) \right]
\]
which, considering the signs of the left hand side and of $g'(\rho_K)-g'(\rho_L)$, in turn implies $\bar \rho_\edge > \rho_K$.
\end{proof}
%


\subsection{Existence of a solution}

The aim of this section is to prove the existence of a solution to the discrete problem under study.
It follows from a topological degree argument, linking by a homotopy the problem at hand to a linear system.

\bigskip
This section begins by a lemma which is used further to prove the strict positivity of the pressure.
%
%
\begin{lmm}\label{max}
Let us consider the following problem:
\begin{equation}
\forall K \in \mesh, \qquad
|K|\,\frac{\varphi_1(p_K)-\varphi_{1}(p^\ast_K)}{\dt}+\sum_{\edge=K|L}\,\vsp\,\varphi_2(p_K)-\vsm\,\varphi_2(p_L)=0
\label{eqk}
\end{equation}
where $\varphi_1$ is an increasing function and $\varphi_2$ is a non-decreasing and non-negative function. 
Suppose that there exists $\bar p$ such that:
\begin{equation}
\varphi_{1}(\bar p) 
+ \dt \, \varphi_2(\bar p) \,\max_{K\in\mesh} \left[ 0,\frac{1}{|K|} \sum_{\edge=K|L} \vs \right]
=\min_{K\in\mesh}[\varphi_1(p^\ast_{K})]
\label{eqbarp}
\end{equation}
Then, $\forall K \in \mesh$, $p_K$ satisfies $ p_K \geq \bar{p}$. 
\end{lmm}

\begin{proof}
Let us assume that there exists a cell $\bar K$ such that $p_{\bar K}=\min_{K \in \mesh} p_K <\bar p$.
Multiplying by $\dt/|\bar K|$ the relation \eqref{eqk} written for $K=\bar K$, we get:
\begin{equation}
\varphi_1(p_{\bar K })
+\frac{\dt}{|\bar K|}\,\sum_{\sigma=\bar K |L} \,(\vsp\,\varphi_2(p_{\bar K})-\vsm\,\varphi_2(p_L))
= \varphi_1(p^\ast_{\bar K})
\label{eqkbar}
\end{equation}
Then substracting from \eqref{eqbarp} we have:
\[
\begin{array}{ll} \displaystyle
\varphi_1(p_{\bar K})-\varphi_1(\bar p)
+\frac{\dt}{|\bar K|} \, \sum_{\edge=\bar K |L} \vsp \,\varphi_2 (p_{\bar K})-\vsm\,\varphi_2(p_L)
& \\ \displaystyle \hfill
- \dt\,\varphi_2 (\bar p)\,\max_{K \in \mesh}\left[0,\frac{1}{|K|} \sum_{\edge=K|L} \vs \right]
& \displaystyle
=\varphi_1(p^\ast_{\bar K})-\min_{K \in \mesh}[\varphi_1(p^\ast_K)]\geq 0
\end{array}
\]
The previous relation can be written as $T_1 + T_2 + T_3 \geq 0$ with:
\[
\begin{array}{l}
T_1= \varphi_1(p_{\bar K})-\varphi_1(\bar p)
\\[2ex] \displaystyle
T_2= \dt \, \varphi_2(p_{\bar K})\left[ \sum_{\edge=\bar K|L} \frac{1}{|\bar K|} \vs \right]
- \dt \,\varphi_{2}(\bar p)\,\max_{K \in \mesh}\left[0,\sum_{\edge=K|L}\frac{1}{|K|} \vs \right]
\\[4ex] \displaystyle
T_3= \frac{\dt}{|\bar K|}\,\sum_{\edge=\bar K |L} \vsm \,(\varphi_2(p_{\bar K})-\varphi_{2}(p_L))
\end{array}
\]
As $\varphi_1(\cdot)$ is an increasing function and, by assumption, $p_{\bar K}<\bar p$, we have $T_1 <0$.
Similarly, $0 \leq \varphi_2(p_{\bar K}) \leq \varphi_2(\bar p)$ and the discrete divergence over $\bar K$ (\ie\ $1/|K|\ \sum_{\edge=\bar K|L} u_\edge$) is necessarily smaller than the maximum of this quantity over the cells of the mesh, thus $T_2 <0$.
Finally, as, by assumption, $p_{\bar K} \leq p_L$ for any neighbouring cell $L$ of $\bar K$, $\varphi_2(\cdot)$ is a non-decreasing function and $\vsm \geq 0$, $T_3 \leq 0$.
We thus obtain a contradiction with the initial hypothesis, which proves $p_K \geq \bar p,\ \forall K \in \mesh$.
\end{proof}
We now state the abstract theorem which will be used hereafter; this result follows from standard arguments of the topological degree theory (see \cite{deimling} for an exposition of the theory and \cite{eym-98-an} for another utilisation for the same objective as here, namely the proof of existence of a solution to a numerical scheme).
\begin{thrm}[A result from the topological degree theory]\label{degre} 
Let $N$ and $M$ be two positive integers and $V$ be defined as follows:
\[
V=\{ (x,y) \in \xR^N \times \xR^M \mbox{ such that } y>0\}
\]
where, for any real number $c$, the notation $y>c$ means that each component of $y$ is greater than $c$.
Let $b\in \xR^N \times \xR^M$ and $f$ and $F$ be two continuous functions respectively from $V$ and $V\times[0,1]$ to $\xR^N \times \xR^M$ satisfying:
\begin{enumerate}
\item[(i)]$F(\cdot,\,1)=f(\cdot)$;
\item[(ii)] $\forall \alpha \in [0,1]$, if $v$ is such that $F(v,\alpha)=b$ then $v \in W$, where $W$ is defined as follows:
\[
W=\{ (x,y) \in \xR^N \times \xR^M \mbox{ s.t. } \Vert x \Vert < C_1 \mbox{ and } \epsilon< y < C_2\}
\]
with $C_1$, $C_2$ and $\epsilon$ three positive constants and $\Vert \cdot \Vert$ is a norm defined over $\xR^N$;
\item[(iii)] the topological degree of $F(\cdot,0)$ with respect to $b$ and $W$ is equal to $d_0 \neq 0$.
\end{enumerate}
Then the topological degree of $F(\cdot,1)$ with respect to $b$ and $W$ is also equal to $d_0 \neq 0$; consequently, there exists at least a solution $v\in W$ such that $f(v)=b$.
\end{thrm}

\bigskip
We are now in position to prove the existence of a solution to the considered discrete Darcy system, for fairly general equations of states.
%
%
\begin{thrm}[Existence of a solution]
Let us suppose that the equation of state $\varrho(\cdot)$ is such that:
\begin{enumerate}
\item $\varrho(\cdot)$ is increasing, $\varrho(0)=0$ and $\displaystyle \lim_{z \rightarrow +\infty} \varrho(z) = +\infty$,
\item there exists an elastic potential $P$ such that (\ref{elasticpot_gene}) holds, the function $z \mapsto z\, P(z)$ is once continuously differentiable and strictly convex and $z\,P(z)\geq -C_P,\ \forall z \in (0,+\infty)$, where $C_P$ is a non-negative constant.
\end{enumerate}
Then, there exists a solution $(u_{\sigma},p_{K})_{\sigma\in\edgesint,\, K \in \mesh}$ to the discrete problem at hand \eqref{pbdisc-disc}.
\label{existence_gene}\end{thrm}

\begin{proof}
This proof is obtained by applying theorem \ref{degre}.
Let $N=d\ {\rm card}(\edgesint)$, $M=\ {\rm card}(\mesh)$ and the mapping $F:\ V \times [0,1] \rightarrow \xR^N \times \xR^M$ be given by:
\[
F(u,p,\alpha)= 
\left| \begin{array}{ll} \displaystyle v_{\edge,i},\ & \edge \in \edgesint,\ 1 \leq i \leq d \\[1ex] q_K,\ & K \in \mesh \end{array} \right.
\]
with:
\begin{equation}
\left| \begin{array}{l} \displaystyle 
v_{\edge,i}=
a(u,\varphi_\edge^{(i)}) - \alpha \int_{\Omega,h} p\ \nabla \cdot \varphi_\edge^{(i)} \dx - \int_\Omega f \cdot \varphi_\edge^{(i)} \dx
\\[2ex] \displaystyle
q_K=\frac{|K|}{\dt} \left[\varrho(p_K)-\varrho(p_K^\ast) \right]
+ \alpha \sum_{\edge=K|L} \vsp\, \varrho_\alpha(p_K) - \vsm\, \varrho(p_L)
\end{array} \right .
\label{eqalpha}
\end{equation}
where, $\forall K \in \mesh$, $p_K^\ast$ is chosen such that $\rho_K^\ast=\varrho(p_K^\ast)$; note that, by assumption, $\varrho(\cdot)$ is one to one from $(0,+\infty)$ to $(0,+\infty)$, so the preceding definition make sense.
The problem $F(u,p,1)=0$ is exactly the system \eqref{pbdisc-disc}.

\medskip
The present proof is built as follows: by applying theorem \ref{VF2}, we obtain a control on $u$ in the discrete norm, uniform with respect to $\alpha$.
Since we work on a finite dimensional space, we then obtain a control on $p$ in $\xLinfty$ by using the conservativity of the system of equations.
For the same reason, the control on $u$ yields a bound in $\xLinfty$ of the value of the discrete divergence, which is shown to allow, by lemma \ref{max}, to bound $p$ away from zero independently of $\alpha$.
The proof finally ends by examining the properties of the system $F(u,p,0)=0$.

\bigskip
\noindent \underline{\textit{Step 1}}: $\alpha\in (0,1]$, $\norms{\cdot}$ estimate for the velocity.
\\[1ex]
Setting $v_\edge=0$ in \eqref{eqalpha}, multiplying the corresponding equation by $u_{\edge,i}$ and summing over $\edge \in \edgesint$ and $1 \leq i\leq d$ yields the following equation:
\[
a(u,u) - \alpha \int_{\Omega,h} p \, \nabla \cdot u \dx = \int_\Omega f \cdot u \dx
\]
By a computation very similar to the proof of theorem \ref{VF2}, we see that, from the second relation of \eqref{eqalpha} with $q_K=0$:
\[
-\alpha \int_{\Omega,h} p \, \nabla \cdot u \dx \geq 
\frac 1 \dt \sum_{K\in\mesh} |K|\ \left[\rho_K\,P(\rho_K) - \rho^\ast_K\,P(\rho\ast_K)\right]
\]
where $\rho_K=\varrho(p_k)$.
By the stability of the bilinear form $a(\cdot,\cdot)$ and Young's inequality, we thus get:
\begin{equation}
\underbrace{\frac{c_{\rm a}} 2 \normsd{u}}_{\displaystyle T_1}
+ \underbrace{\frac 1 \dt \sum_{K\in\mesh} |K|\ \rho_K\,P(\rho_K)}_{\displaystyle T_2}
\leq
\frac 1 {2 c_{\rm a}}\normSd{f}
+ \frac 1 \dt \sum_{K\in\mesh} |K|\ \rho^\ast_K\,P(\rho\ast_K)
\label{e_est}\end{equation}
By assumption, $T_2 \geq -C_p |\Omega|$ and we thus get the following estimate on the discrete norm of the velocity:
\begin{equation}
\norms{u}\leq C_{1}
\label{est_u}\end{equation}
where $C_1$ only depends on the data of the problem, \ie\ $\dt$, $f$ and $\rho^\ast$ and not on $\alpha$.

\bigskip
\noindent \underline{\textit{Step 2}}: $\alpha\in (0,1]$, $\xLinfty$ estimate for the pressure.
\\[1ex]
Let us now turn to the estimate of the pressure.
By conservativity of the discrete mass balance, it is easily seen that:
\[
\sum_{K\in \mesh} |K|\ \varrho(p_K)=\sum_{K\in \mesh} |K|\ \rho^\ast_K
\]
As each term in the sum on the left hand side is non-negative, we thus have:
\[
\forall K \in \mesh, \qquad \varrho(p_K) \leq \frac 1 {\min_{K\in\mesh} |K|}\ \sum_{K\in \mesh} |K|\ \rho^\ast_K
\]
which, as, by assumption, $\displaystyle \lim_{z \rightarrow +\infty} \varrho(z) = +\infty$, yields:
\begin{equation}
\forall K \in \mesh, \qquad p_K \leq C_2
\label{est_p_max}\end{equation}
where $C_2$ only depends on the data of the problem.

\bigskip
\noindent \underline{\textit{Step 3}}: $\alpha\in (0,1]$, $p$ bounded away from zero.
\\[1ex]
Applying lemma \ref{max} with $\varphi_1(\cdot)=\varphi_2(\cdot)=\varrho(\cdot)$, we get:
\[
\forall K \in \mesh, \qquad p_K \geq \bar p_\alpha
\]
where $\bar p_\alpha$ is given by:
\[
\varrho(\bar p_\alpha)=\frac{\displaystyle \min_{K\in\mesh} \varrho(p^\ast)}
{\displaystyle 1+\dt \max_{K \in \mesh}\left[0,\alpha \sum_{\edge=K|L}\frac{1}{|K|}\,\vs \right]}
\]
Note that $\bar p_\alpha$ is well defined since, by assumption, $\varrho(\cdot)$ is one to one from $(0,+\infty)$ to $(0,+\infty)$.
As $\alpha \leq 1$, we get:
\[
\varrho(\bar p_\alpha) \geq \frac{\displaystyle \min_{K\in\mesh} \varrho(p^\ast)}
{\displaystyle 1+\dt \max_{K \in \mesh}\left[0,\sum_{\edge=K|L}\frac{1}{|K|}\,\vs \right]}
\]
and, by equivalence of norms in a finite dimensional space, the bound \eqref{est_u} also yields a bound in the $\xLinfty$ norm and, finally, an upper bound for the right hand side of this relation.
We thus get, still since $\varrho(\cdot)$ is increasing on $(0,+\infty)$ that, $\forall \alpha \in (0,1]$, $\bar p_\alpha \geq \epsilon_1$, and, finally:
\begin{equation}
\forall K \in \mesh, \qquad p_K \geq \epsilon_1
\label{est_p_min}\end{equation}
where $\epsilon$ only depends on the data.

\bigskip
\noindent \underline{\textit{Step 4}}: conclusion.
\\[1ex]
For $\alpha=0$, the system $F(u,p,0)=0$ reads:
\[
\left| \begin{array}{ll} \displaystyle 
a(u,\varphi_\edge^{(i)}) = \int_\Omega f \cdot \varphi_\edge^{(i)} \dx \hspace{10ex}
&
\forall \edge \in \edgesint,\ 1\leq i \leq d
\\[3ex] \displaystyle
\varrho(p_K)=\varrho(p_K^\ast)
&
\forall K \in \mesh
\end{array} \right .
\]
Since $\varrho(\cdot)$ is one to one from $(0,+\infty)$ to $(0,+\infty)$ and thanks to the stability of the bilinear form $a(\cdot,\cdot)$, this system has one and only one solution (which, for the pressure, reads of course $p_K=p_K^\ast,\ \forall K \in \mesh$), which satisfies:
\begin{equation}
\norms{u}\leq C_{3},\qquad
\epsilon_2=\min_{K\in\mesh} p_K^\ast \leq p \leq \max_{K\in\mesh} p_K^\ast=C_4
\label{est_up}\end{equation}
In addition, the jacobian of this system is block diagonal: the first block, associated to the first relation, is constant (this part of the system is linear) and non-singular; the second one, associated to the second relation, is diagonal, and each diagonal entry is equal to the derivative of $\varrho(\cdot)$, taken at the considered point.
As the function $\varrho(\cdot)$ is increasing, this jacobian matrix does not vanish for the solution of the system.
Let $W$ be defined by:
\[
W=\{(u,p) \in \xR^N \times \xR^M \mbox{ such that } \norms{\tilde\rho}{u} < 2 \max(C_1,C_3)
\mbox{ and } \frac 1 2\,\min(\epsilon_1,\epsilon_2)< p < 2 \max(C_1,C_3) \}
\]
The topological degree of $F(\cdot,\cdot,0)$ with respect to $0$ and $W$ does not vanish, and by inequalities \eqref{est_u}, \eqref{est_p_min}, \eqref{est_p_max} and \eqref{est_up}, theorem \ref{degre} applies, which concludes the proof.
\end{proof}


\subsection{Some cases of application}

First of all, let us give some examples for the bilinear form $a(\cdot,\cdot)$, for which the theory developped in this work holds.
The first of them is:
\[
a(u,v)=\int_\Omega u \cdot v \dx, \qquad \norms{u}=\normLdv{u}, \qquad \normS{f}=\normLdv{f}
\]
This choice for $a(\cdot,\cdot)$ yields a discrete Darcy-like problem which is, up to numerical integration technicalities, the projection step arising in the pressure correction scheme which is considered in the present paper (see section \ref{sec:scheme}).
Note that, in this case, the boundary condition $u\in\xHone_0(\Omega)^d$ does not make sense at the continuous level; in addition, the considered discretization is known to be not consistent enough to yield convergence (see remark \ref{non-consistent!} hereafter) for the Darcy problem.

\medskip
The bilinear form associated to the Stokes problem provides another example of application.
It may read in this case:
\[
a(u,v)=\int_{\Omega,h} \nabla u \cdot \nabla v \dx, \qquad \norms{u}=|u|_{\xHone(\Omega)^d}
\]
or, without additional theoretical difficulties:
\[
a(u,v)= \mu \int_{\Omega,h} \nabla u \cdot \nabla v \dx + \frac \mu 3 \int_{\Omega,h} (\nabla \cdot u) \ (\nabla \cdot v) \dx
\] 
this latter form, where the real number $\mu >0$ is the viscosity, corresponding to the physical shear stress tensor expression for a compressible flow of a constant viscosity Newtonian fluid.

\medskip
In addition, consider the steady Navier-Stokes equations, or, more generally, a time step of a (semi-)implicit time discretization of the unsteady Navier-Stokes equations, in which case $a(\cdot,\cdot)$ and $f$ reads:
\[
\begin{array}{l} \displaystyle
a(u,v) = \frac 1 \dt \int_\Omega \rho u \cdot v \dx + \int_{\Omega,h} (\nabla \cdot \rho w \otimes u) \cdot v \dx
+ \mu \int_{\Omega,h} \nabla u \cdot \nabla v \dx + \frac \mu 3 \int_{\Omega,h} (\nabla \cdot u) \ (\nabla \cdot v) \dx
\\[3ex] \displaystyle
f= \frac 1 \dt \rho^\ast u^\ast + f_0
\end{array}
\]
where the steady case is obtained with $\dt=+\infty$, $f_0$ is the physical forcing term, $\rho^\ast$ and $u^\ast$ stands for known density and velocity fields and $w$ is an advection field, which may be $u$ itself (and must be $u$ in the steady case) or be derived from the velocity obtained at the previous time steps.
Let us suppose that the following identity holds:
\[
\frac 1 \dt \int_\Omega (\rho u - \rho^\ast u^\ast) \cdot u \dx + \int_{\Omega,h} (\nabla \cdot \rho w \otimes u) \cdot u \dx
\geq \frac 1 {2 \dt} \left[ \int_\Omega \rho |u|^2 - \int_\Omega \rho^\ast |u^\ast|^2 \right]
\]
which is the discrete counterpart of equation \eqref{stab-arg}-$(i)$.
The algorithm considered in this paper provides an example where this condition is verified (see section \ref{sec:scheme}).
Then the present theory applies with little modifications: in the proof of existence theorem \ref{existence_gene}, the right hand side of the preceding equation must be multiplied by the homotopy parameter $\alpha$ (an thus this term vanishes at $\alpha=0$, which yields the problem  considered in step $4$ above); the (uniform with respect to $\alpha$) stability in step 1 stems from the diffusion term, and steps $2$ and $3$ remain unchanged.

\medskip
Note that, in the steady state case, an additional constraint is needed for the problem to be well posed, namely to impose the total mass $M$ of fluid in the computational domain to a given value.
This constraint can be simply enforced by solving an approximate mass balance which reads:
\[
c(h) \, \left[ \rho - \frac M {|\Omega|} \right] + \nabla \cdot \rho u =0
\]
where $|\Omega|$ stands for the measure of $\Omega$, $h$ is the spatial discretization step and $c(h)>0$ must tend to zero with $h$, fast enough to avoid any loss of consistency.
With this form of the mass balance, the theory developped here directly applies.

\bigskip
Examining now the assumptions for the equation of state in theorem \ref{existence_gene}, we see that our results hold with equations of state of the form:
\[
\varrho(p)= p^{1/\gamma} \qquad \mbox{or, equivalently} \qquad \rho=p^\gamma, \qquad \mbox{where } \gamma > 1
\]
In this case, the elastic potential is given by equation \eqref{elasticpot}, which yields:
\[
P(\rho)=\frac 1 {\gamma -1}\ \rho^{\gamma -1},\qquad \rho P(\rho)= \frac 1 {\gamma -1}\ \rho^\gamma \quad (=\frac 1 {\gamma -1}\ p)
\]
The same conclusion still holds with $\gamma=1$ (\ie \ $p=\rho$), with $P(\rho)=\log(\rho)$ satisfying equation \eqref{elasticpot_gene}.
The case $\gamma > 1$ is for instance encountered for isentropic perfect gas flows, whereas $\gamma=1$ corresponds to the isothermal case.
It is worth noting that this range of application is larger than what is known for the continuous case, for which the existence of a solution is known only in the case $\gamma > d/2$ \cite{pll-98-mat, fei-04-dyn, nov-04-int}.


\section{A pressure correction scheme} \label{sec:scheme}

In this section, we build a pressure correction numerical scheme for the solution of compressible barotropic Navier-Stokes equations \eqref{nsb}, based on the low order non-conforming finite element spaces used in the previous section, namely the Crouzeix-Raviart or Rannacher-Turek elements.

\medskip
The presentation is organized as follows.
First, we write the scheme in the time semi-discrete setting (section \ref{subsec:sd}).
Then we prove a general stability estimate which applies to the discretization by a finite volume technique of the convection operator (section \ref{subsec:conv}).
The proposed scheme is built in such a way that the assumptions of this stability result hold (section \ref{subsec:mom}); this implies first to perform a prediction of the density, as a non-standard first step of the algorithm and, second, to discretize the convection terms in the momentum balance equation by a finite volume technique specially designed to this purpose.
The discretization of the projection step (section \ref{subsec:proj}) also combines the finite element and finite volume methods, in such a way that the theory developped in section \ref{sec:darcy} applies; in particular, the proposed discretization allows to take benefit of the pressure or density control induced by the pressure work, \ie\ to apply theorem \ref{VF2}.
The remaining steps of the algorithm are described in section \ref{subsec:renorm} and an overview of the scheme is given in section \ref{subsec:overview}.
The following section (section \ref{subsec:stab}) is devoted to the proof of the stability of the algorithm.
Finally, we shortly address some implementation difficulties (section \ref{subsec:impl}), then provide some numerical tests (section \ref{subsec:num-exp}) performed to assess the time and space convergence of the scheme.

\subsection{Time semi-discrete formulation}\label{subsec:sd}

Let us consider a partition $0=t_0 < t_1 <\ldots < t_n=T$ of the time interval $(0,T)$, which, for the sake of simplicity, we suppose uniform.
Let $\dt$ be the constant time step $\dt=t_{k+1}-t_k$ for $k=0,1,\ldots,n$. 
In a time semi-discrete setting, the scheme considered in this paper reads:
\vskip-.6cm
\begin{eqnarray}
&& \mbox{1 -- Solve for } \tilde \rho^{n+1}\mbox{\ : \quad}
\frac{\tilde \rho^{n+1}-\rho^{n}}{\dt}+\nabla \cdot (\tilde \rho^{n+1}\,u^{n})=0
\label{stab-scheme-1}\\
&& \mbox{2 -- Solve for } \tilde p^{n+1}\mbox{\ : \quad}
-\nabla \cdot \left(\frac{1}{\tilde \rho ^{n+1}} \nabla \tilde p^{n+1}\right) = 
-\nabla \cdot \left(\frac{1}{\sqrt{\tilde \rho^{n+1}}\sqrt{\tilde \rho^{n}}}\nabla p^n \right)
\label{stab-scheme-2}\\
&& \mbox{3 -- Solve for } \tilde u^{n+1}\mbox{\ : \quad}
\nonumber \\
&& \quad
\frac{\tilde \rho^{n+1}\, \tilde u^{n+1}-\rho^{n}\,u^{n}}{\dt}
+\nabla \cdot (\tilde \rho^{n+1}\,u^n\otimes\tilde u^{n+1}) 
+\nabla \tilde p^{n+1}
-\nabla \cdot \tau(\tilde u^{n+1}) 
= f^{n+1}
\label{stab-scheme-3}\\
&& \mbox{4 -- Solve for } \bar u^{n+1},\, p^{n+1},\, \rho^{n+1}\mbox{\ : \quad}
\nonumber \\
&& \qquad \left| \begin{array}{l} \displaystyle
\tilde \rho^{n+1}\,\frac{\bar u^{n+1}- \tilde u^{n+1}}{\dt}
+ \nabla (p^{n+1}-\tilde p^{n+1})=0
\\ \displaystyle
\frac{\varrho(p^{n+1})-\rho^n}{\dt}
+ \nabla \cdot \left( \varrho(p^{n+1})\, \bar u^{n+1} \right)=0
\\ \displaystyle
\rho^{n+1}=\varrho(p^{n+1})
\end{array} \right .
\label{stab-scheme-4}\\
&& \mbox{5 -- Compute } u^{n+1} \mbox{ given by: \quad}
\sqrt{\rho^{n+1}}\,u^{n+1}=\sqrt{\tilde \rho^{n+1}}\bar u^{n+1}
\label{stab-scheme-5}
\end{eqnarray}

The first step is a prediction of the density, used for the discretization of the time derivative of the momentum.
As remarked by Wesseling \textit{et al} \cite{bij-98-uni, wes-01-pri}, this step can be avoided when solving the Euler equations: in this case,  the mass flowrate may be chosen  as an unknown, using the explicit velocity as an advective field in the discretization of the convection term in the momentum balance;   the velocity is then updated  by dividing by the density at the end of the time step.
For viscous flows, if the discretization of the diffusion term is chosen to be implicit, both the mass flowrate and the velocity appear as unknowns in the momentum balance; this seems to impeed the use of this trick.
Let us emphasize that the way the step one is carried out (\textit{i.e.} solving a discretization of the mass balance instead as, for instance, performing a Richardson's extrapolation) is crucial for the stability.

\medskip
Likewise, the second step is a renormalization of the pressure the interest of which is clarified only by the stability analysis.
A similar technique has already been introduced by Guermond and Quartapelle for variable density incompressible flows \cite{gue-00-proj}.

\medskip
Step 3 consists in a classical semi-implicit solution of the momentum equation to obtain a predicted velocity.

\medskip
Step 4 is a nonlinear pressure correction step, which degenerates in the usual projection step as used in incompressible flow solvers when the density is constant (e.g. \cite{mar-98-nav}).
Taking the divergence of the first relation of (\ref{stab-scheme-4}) and using the second one to eliminate the unknown velocity $\bar u^{n+1}$ yields a non-linear elliptic problem for the pressure.
This computation is formal in the semi-discrete formulation, but, of course, is necessarily made clear at the algebraic level, as described in section \ref{subsec:impl}.
Once the pressure is computed, the first relation yields the updated velocity and the third one gives the end-of-step density.


\subsection{Stability of the advection operator: a finite-volume result} \label{subsec:conv}

The aim of this section is to state and prove a discrete analogue to the stability identity \eqref{stab-arg}-$(i)$, which may be written for any sufficiently regular functions $\rho$, $z$ and $u$ as follows:
\[
\int_\Omega \left[ \frac{\partial \rho z}{\partial t} + \nabla \cdot (\rho z u)\right ] z \dx
= \frac{1}{2} \, \frac{d}{dt} \int_\Omega \rho z^2 \dx
\label{VF1_continu}
\]
and holds provided that the following balance is satisfied by $\rho$ and $u$:
\[
\frac{\partial \rho}{\partial t} + \nabla \cdot (\rho u) =0
\]
As stated in introduction, applying this identity to each component of the velocity yields the central argument of the proof of the kinetic energy theorem.

\bigskip
The discrete analogue to this identity is the following.
\begin{thrm}[Stability of the advection operator]\label{VF1}Let $(\rho_K^\ast)_{K\in \mesh}$ and $(\rho_K)_{K\in \mesh}$ be two families of positive real number satisfying the following set of equation:
\begin{equation}
\forall K \in \mesh,\qquad \frac{|K|}{\dt} \ (\rho_K - \rho^\ast_K) + \sum_{\edge=K|L} \fluxK=0
\label{mass_bal}\end{equation}
where $\fluxK$ is a quantity associated to the edge $\edge$ and to the control volume $K$; we suppose that, for any internal edge $\edge=K|L$, $\fluxK=-\fluxL$.
Let $(z_K^\ast)_{K\in \mesh}$ and $(z_K)_{K\in \mesh}$ be two families of real numbers.
For any internal edge $\edge=K|L$, we define $z_\edge$ either by $z_\edge=\frac 1 2 (z_K +z_L)$, either by $z_\edge=z_K$ if $\fluxK \geq 0$ and $z_\edge=z_L$ otherwise.
The first choice is referred to as an the "centered choice", the second one as "the upwind choice".
In both cases, the following stability property holds:
\begin{equation}
\sum_{K \in \mesh} z_K \left[ \frac{|K|}\dt\ (\rho_K\, z_K -\rho_K^\ast\, z_K^\ast)+
\sum_{\edge=K|L} \fluxK \, z_\edge \right]
\geq
\frac{1}{2}\ \sum_{K\in\mesh} \frac{|K|}\dt\ \left[ \rho_K\, z_K^2 -\rho_K^\ast\, {z_K^\ast}^2\right]
\label{ecinetique} \end{equation}
\end{thrm}

\begin{proof}
We write:
\[
\sum_{K \in \mesh} z_K \left[ \frac{|K|}\dt\ (\rho_K\, z_K -\rho_K^\ast\, z_K^\ast)+
\sum_{\edge=K|L} \fluxK \, z_\edge \right]
=T_1 + T_2
\]
where $T_1$ and $T_2$ reads:
\[
T_1= \sum_{K \in \mesh} \frac{|K|}\dt\ z_K\ (\rho_K\, z_K -\rho_K^\ast\, z_K^\ast), \hspace{10ex}
T_2=\sum_{K \in \mesh} z_K \left[ \sum_{\edge=K|L} \fluxK \, z_\edge \right]
\]
The first term reads:
\[
T_1= \sum_{K \in \mesh} \frac{|K|}\dt\ \left[z_K^2\, (\rho_K-\rho_K^\ast) + \rho_K^\ast \, z_K\, (z_K -z_K^\ast) \right]
\]
Developping the last term by the identity $a(a-b)=\frac 1 2 ( a^2 + (a-b)^2 -b^2)$, we get:
\[
T_1= \underbrace{\sum_{K \in \mesh} \frac{|K|}\dt\ z_K^2\, (\rho_K-\rho_K^\ast)}_{\displaystyle T_{1,1}}
+ \underbrace{\frac 1 2 \sum_{K \in \mesh} \frac{|K|}\dt\ \rho_K^\ast\, (z_K^2 - {z_K^\ast}^2)}_{\displaystyle T_{1,2}}
+ \underbrace{\frac 1 2 \sum_{K \in \mesh} \frac{|K|}\dt\ \rho_K^\ast\, (z_K - z_K^\ast)^2}_{\displaystyle T_{1,3}}
\]
The last term, namely $T_{1,3}$, is always non-negative and can be seen as a dissipation associated to the backward time discretization of equation \eqref{ecinetique}.
We now turn to $T_2$:
\[
T_2=\underbrace{\sum_{K \in \mesh} z_K^2 \left[ \sum_{\edge=K|L} \fluxK\right]}_{\displaystyle T_{2,1}}+
    \underbrace{\sum_{K \in \mesh} z_K \left[ \sum_{\edge=K|L} \fluxK \, (z_\edge-z_K) \right]}_{\displaystyle T_{2,2}}
\]
The first term, namely $T_{2,1}$, will cancel with $T_{1,1}$ by equation \eqref{mass_bal}.
The second term reads, developping as previously the quantity $z_K\, (z_\edge-z_K)$:
\[
T_{2,2}= -\frac 1 2 \sum_{K \in \mesh} z_K^2 \left[ \sum_{\edge=K|L} \fluxK \right]
- \underbrace{\frac 1 2 \sum_{K \in \mesh} \left[ \sum_{\edge=K|L} \fluxK\ [(z_\edge-z_K)^2 - z_\edge^2] \right]}_{\displaystyle T_{2,3}}
\]
Reordering the sum in the last term, we have, as $\fluxK=-\fluxL$:
\[
T_{2,3}=\frac 1 2 \sum_{\edge \in \edgesint\ (\edge=K|L)} \fluxK\ [(z_\edge-z_K)^2-(z_\edge-z_L)^2]
\]
This expression can easily be seen to vanish with the centered choice.
With the upwind choice, supposing without loss of generality that we have chosen for the edge $\edge=K|L$ the orientation such that $\fluxK \geq 0$, we get, as $z_\edge=z_K$:
\[
T_{2,3}=- \frac 1 2 \sum_{\edge \in \edgesint\ (\edge=K|L)} \fluxK\ (z_K-z_L)^2 \leq 0
\]
We thus have, by equation \eqref{mass_bal}:
\[
T_{2,2} \geq -\frac 1 2 \sum_{K \in \mesh} z_K^2 \left[ \sum_{\edge=K|L} \fluxK \right]
=\frac 1 2 \sum_{K \in \mesh} \frac{|K|}\dt z_K^2\, (\rho_K - \rho_K^\ast)
\]
and thus:
\[
T_1 + T_2 \geq \frac 1 2 \sum_{K \in \mesh} \frac{|K|}\dt \left[z_K^2\, (\rho_K - \rho_K^\ast) + \rho_K^\ast\, (z_K^2 - {z_K^\ast}^2)\right]
\]
which concludes the proof.
\end{proof}

\begin{rmrk}
Equation \eqref{mass_bal} can be seen as a discrete mass balance, with $\fluxK$ standing for the mass flux across the edge $\edge$, and the right hand side of \eqref{ecinetique} may be derived by the multiplication by $z_K$ and summation over the control volumes of the transport terms in a discrete balance equation for the quantity $\rho z$, reading:
\[
\forall K \in \mesh, \qquad 
\frac{|K|}\dt\ (\rho_K\, z_K -\rho_K^\ast\, z_K^\ast)+ \sum_{\edge=K|L} \fluxK \, z_\edge +\cdots \mbox{[possible diffusion terms]} \cdots =0
\]
In this context, the relation \eqref{mass_bal} is known to be exactly the compatibility condition which ensures a discrete maximum principle for the solution $z$ of this transport equation, provided that the upwind choice (or any monotone choice) is made for the expression of $z_\edge$ \cite{lar-91-how}.
We proved here that the same compatibility condition ensures a $\xLtwo$ stability for $\rho^{1/2} z$. 
\end{rmrk}


\subsection{Spatial discretization of the density prediction and the momentum balance equation}\label{subsec:mom}

The main difficulty in the discretization of the momentum balance equation is to build a discrete convection operator which enjoys the stability property \eqref{stab-arg}-$(i)$.
To this purpose, we derive for this term a finite volume discretization which satisfies the assumptions of theorem \ref{VF1}.

\medskip
The natural space discretization for the density is the same as for the pressure, \ie \ piecewise constant functions over each element.
For the Rannacher-Turek element, this legitimates a standard mass lumping technique for the time derivative term, since no additional accuracy seems to have to be expected from a more precise integration.
For the Crouzeix-Raviart element, the mass matrix is already diagonal.
Let the quantity $|D_\edge|$ be defined as follows:
\begin{equation}
|D_\edge| \eqdef \int_\Omega \varphi_\edge \dx >0
\label{int-phis}\end{equation}
For the Crouzeix-Raviart element, $|D_\edge|$ can be identified to the measure of the cone with basis $\edge$ and with the mass center of the mesh as opposite vertex.
The same property holds for the Rannacher-Turek element in the case of quandrangles ($d=2$) or cuboids ($d=3$), which are the only case considered here, even though extensions to non-perpendicular grids are probably possible.

\begin{figure}
\setlength{\unitlength}{1cm}
\begin{center}
\begin{picture}(12.8,6)
\thicklines
\put(1.5,0){\line(1,0){7}}
\put(1.5,3){\line(1,0){7}}
\put(1.5,6){\line(1,0){7}}
\put(2,-0.2){\line(0,1){6.4}}
\put(8,-0.2){\line(0,1){6.4}}

\thinlines
\put(3.,0.5){\line(2,1){5}}
\put(8,0){\line(-2,1){6}}
\put(2,3){\line(2,1){6}}
\put(8,3){\line(-2,1){5}}
\thicklines
\put(5,1.5){\line(2,1){3}}
\put(5,1.5){\line(-2,1){3}}
\put(2,3){\line(2,1){3}}
\put(8,3){\line(-2,1){3}}
\thinlines
\put(2.2,2.9){\line(0,1){0.2}}
\put(2.4,2.8){\line(0,1){0.4}}
\put(2.6,2.7){\line(0,1){0.6}}
\put(2.8,2.6){\line(0,1){0.8}}
\put(3.0,2.5){\line(0,1){1.0}}
\put(3.2,2.4){\line(0,1){1.2}}
\put(3.4,2.3){\line(0,1){1.4}}
\put(3.6,2.2){\line(0,1){1.6}}
\put(3.8,2.1){\line(0,1){1.8}}
\put(4.0,2.0){\line(0,1){2.0}}
\put(4.2,1.9){\line(0,1){2.2}}
\put(4.4,1.8){\line(0,1){2.4}}
\put(4.6,1.7){\line(0,1){2.6}}
\put(4.8,1.6){\line(0,1){2.8}}
\put(5.0,1.5){\line(0,1){3.0}}
\put(5.2,1.6){\line(0,1){2.8}}
\put(5.4,1.7){\line(0,1){2.6}}
\put(5.6,1.8){\line(0,1){2.4}}
\put(5.8,1.9){\line(0,1){2.2}}
\put(6.0,2.0){\line(0,1){2.0}}
\put(6.2,2.1){\line(0,1){1.8}}
\put(6.4,2.2){\line(0,1){1.6}}
\put(6.6,2.3){\line(0,1){1.4}}
\put(6.8,2.4){\line(0,1){1.2}}
\put(7.0,2.5){\line(0,1){1.0}}
\put(7.2,2.6){\line(0,1){0.8}}
\put(7.4,2.7){\line(0,1){0.6}}
\put(7.6,2.8){\line(0,1){0.4}}
\put(7.8,2.9){\line(0,1){0.2}}

\put(2.2,0.1){control volume K}
\put(2.2,5.65){control volume L}
\put(8.5,4.){edge $\sigma=K|L$}
\put(8.5,5.){diamond $D_\sigma$}
\put(8.4,4.05){\vector(-1,-1){1.05}}
\put(8.4,5.05){\vector(-2,-1){3}}
\end{picture}
\end{center}
\caption{Dual finite volume mesh: the so-called "diamond cells".}
\label{diamonds}
\end{figure}
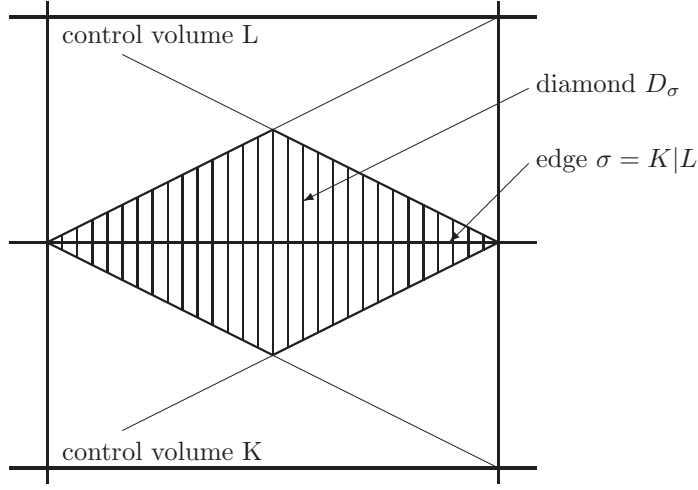

\medskip
For each internal edge $\edge=K|L$, this conic volume is denoted by $D_{K,\edge}$; the volume $D_\edge= D_{K,\edge}\cup D_{L,\edge}$ is referred to as the "diamond cell" associated to $\edge$, and $D_{K,\edge}$ is the half-diamond cell associated to $\edge$ and $K$ (see figure \ref{diamonds}).
The measure of $D_{K,\edge}$ is denoted by $|D_{K,\edge}|$.

\medskip
The discretization of the term $\rho^n \, u^n$ thus leads, in the equations associated to the velocity on $\edge$, to an expression of the form $\rho^n_\edge \, u^n_\edge$, where $\rho^n_\edge$ results from an average of the values taken by the density in the two elements adjacent to $\edge$, weighted by the measure of the half-diamonds:
\begin{equation}
\forall \edge \in \edgesint,\qquad |D_\edge|\ \rho^n_\edge= |D_{K,\edge}|\ \rho^n_K + |D_{L,\edge}|\ \rho^n_L
\label{rho_faces}\end{equation}

\bigskip
In order to satisfy the compatibility condition introduced in the previous section, a prediction of the density is first performed, by a finite volume discretization of the mass balance equation, taking the diamond cells as control volumes:
\begin{equation}
\forall \edge \in \edgesint,\qquad
\frac{|D_\edge|} \dt\ (\tilde \rho_\edge^{n+1}- \rho^n_\edge) + \sum_{\edged \in \edges(D_\edge)} \fluxd^{n+1}=0
\label{mass-diamond}\end{equation}
where $\edges(D_\edge)$ is the set of the edges of $D_\edge$ and $\fluxd$ stands for the mass flux across $\edged$ outward $D_\edge$.
This latter quantity is expressed as follows:
\[
F_{\edged}^{n+1}=|\edged|\ u_{\edged}^n\cdot n_{\edged,\edge}\ \tilde \rho_{\edged}^{n+1}
\]
where $|\edged|$ is the measure of $\edged$, $n_{\edged,\edge}$ is the normal to $\edged$ outward $D_\edge$, the velocity $u_{\edged}^n$ is obtained by interpolation at the center of $\edged$ in $u^n$ (using the standard finite element expansion) and $\tilde \rho_{\edged}^{n+1}$ is the density at the edge, calculated by the standard upwinding technique (\ie\ either $\tilde \rho_\edge^{n+1}$ if $u_{\edged}^n\cdot n_{\edged,\edge} \geq 0$ or $\tilde \rho_{\edge'}^{n+1}$ otherwise, with $\edge'$ such that $\edged$ separates $D_\edge$ and $D_{\edge'}$, which we denote by $\edged=D_\edge |  D_{\edge'}$).

\bigskip
The discretization convection terms of the momentum balance equation are built from relation \eqref{mass-diamond}, according to the structure which is necessary to apply theorem \ref{VF1}.
This yields the following discrete momentum balance equation:
\begin{equation}
\begin{array}{l} \displaystyle
\frac{|D_\edge|} \dt\ (\tilde \rho_\edge^{n+1} \tilde u_{\edge,i}^{n+1}- \rho^n_\edge u_{\edge,i}^n)
+ \sum_{\stackrel{\scriptstyle \edged \in \edges(D_\edge),}{\scriptstyle \edged=D_\edge |  D_{\edge'}}}\frac 1 2\ \fluxd^{n+1}\ (\tilde u_{\edge,i}^{n+1} + u_{\edge',i}^{n+1})
\\ \hspace{8ex}\displaystyle
+ \int_{\Omega,h} \tau(\tilde u^{n+1}) : \nabla \varphi_\edge^{(i)} \dx
-\int_{\Omega,h} \tilde p^{n+1}\ \nabla \cdot \varphi_\edge^{(i)} \dx 
= \int_\Omega f^{n+1} \cdot \varphi_\edge^{(i)}
\hspace{10ex} \edge \in \edgesint,\ 1 \leq i \leq d
\end{array}
\label{momentum}\end{equation}
where we recall that $\varphi_\edge^{(i)}$ the vector shape function associated to $v_{\edge, i}$, which reads $\varphi_\edge^{(i)}=\varphi_\edge \, e_i$ with $e_i$ is the $i^{th}$ vector of the canonical basis of $\xR^d$ and $\varphi_\edge$ the scalar shape function, and the notation $\int_{\Omega,h}$ means $\sum_{K\in\mesh} \int_K$.

\bigskip
Note that, for Crouzeix-Raviart elements, a combined finite volume/finite element method, similar to the technique employed here for the discretization of the momentum balance, has already been analysed for a transient non-linear convection-diffusion equation by Feistauer and co-workers \cite{ang-98-ana, dol-02-err, fei-03-mat}.

\subsection{Spatial discretization of the projection step}\label{subsec:proj}

The fully discrete projection step used in the described algorithm reads:
\begin{equation}
\left| \begin{array}{ll} \displaystyle
|D_\edge|\,\frac{\tilde \rho_\edge^{n+1}}{\dt} \, (\bar u_{\edge,i}^{n+1} - \tilde u_{\edge,i}^{n+1}) 
+ \int_{\Omega,h} (p^{n+1}-\tilde p^{n+1})\ \nabla \cdot \varphi_\edge^{(i)} \dx = 0
\hspace{12ex} & \displaystyle
\edge \in \edgesint,\ 1 \leq i \leq d
\\[2ex] \displaystyle
|K|\ \frac{\varrho(p^{n+1}_K)-\rho^n_K}{\dt} 
+ \sum_{\edge=K|L} (\vsp)^{n+1}\, \varrho(p_K^{n+1}) - (\vsm)^{n+1}\, \varrho(p_L^{n+1})=0
& \displaystyle
K \in \mesh
\end{array}\right .
\label{darcy-disc-un}\end{equation}
where $(\vsp)^{n+1}$ and $(\vsm)^{n+1}$ stands respectively for $\max (\vs^{n+1},\ 0)$ and $\max ( -\vs^{n+1},\ 0)$ with $\vs^{n+1}=|\edge|\, \bar u_\edge^{n+1} \cdot n_{KL}$.
The first (vector) equation may be seen as the finite element discretization of the first relation of the projection step \eqref{stab-scheme-4}, with the same lumping of the mass matrix for the Rannacher-Turek element as in the prediction step.
As the pressure is piecewise constant, the finite element discretization of the second relation of \eqref{stab-scheme-4}, \ie\ the mass balance, is equivalent to a finite volume formulation, in which we introduce the standard first-order upwinding.
Exploiting the expression of the velocity and pressure shape functions, the first set of relations of this system can be alternatively written as follows: 
\begin{equation}
|D_\edge|\,\frac{\tilde \rho_\edge^{n+1}}{\dt} \, (\bar u_\edge^{n+1}-\tilde u_\edge^{n+1}) 
+ |\edge| \, \left[ (p_K^{n+1}-\tilde p_K^{n+1})-(p_L^{n+1}-\tilde p_L^{n+1}) \right]\, n_{KL} = 0
\qquad 
\forall \edge\in\edgesint, \ \edge=K|L
\label{darcy-disc-deux}\end{equation}
or, in an algebraic setting:
\begin{equation}
\frac{1}{\dt} \, \matM_{\tilde \rho^{n+1}} \, (\bar u^{n+1} - \tilde u^{n+1}) + \matB^t \, (p^{n+1}-\tilde p^{n+1}) = 0
\label{darcy-alg}\end{equation}
In this relation, $\matM_{w}$ stands for the diagonal mass matrix  weighted by $(w_\edge)_{\edge \in \edgesint}$ (so, for $1\leq i \leq d$ and $\edge \in \edgesint$, the corresponding entry on the diagonal of $\matM_{\tilde \rho^{n+1}}$ reads  $(\matM_{\tilde \rho^{n+1}})_{\edge,i}= (|D_\edge|\, \tilde \rho_\edge^{n+1})$), $\matB^t$ is the matrix of $\xR^{(d\,N) \times M}$, where $N$ is the number of internal edges (\ie\ $N={\rm card}\,(\edgesint)$) and $M$ is the number of control volumes in the mesh (\ie\ $M={\rm card}\,(\mesh)$), associated to the gradient operator; consequently, the matrix $\matB$ is associated to the opposite of the divergence operator.
Throughout this section, we use the same notation for the discrete function (defined as usual in the finite element context by its expansion using the shape functions) and for the vector gathering the degrees of freedom; so, in relation \eqref{darcy-alg}, $\bar u$ stands for the vector of $\xR^{d\,N}$ components $\bar u_{\edge,i},\ 1\leq i \leq d,\ \edge \in \edgesint$ and $p$ stand for the vector of $\xR^M$ of components $p_K,\ K \in \mesh$.
Both forms \eqref{darcy-disc-deux} and \eqref{darcy-alg} are used hereafter.

\bigskip
We have the following existence result.
\begin{prpstn}
Let the equation of state $\varrho(\cdot)$ be such that $\varrho(0)=0$, $\lim_{z \rightarrow + \infty} \varrho(z)= + \infty$ and there exists an elastic potential function $P(\cdot)$ such that the function $z \mapsto z\,P(z)$ is bounded from above in $(0, + \infty)$,  once continuously differentiable and strictly convex.
Then the nonlinear system \eqref{darcy-disc-un} admits at least one solution and any possible solution is such that $p_K >0$, $\forall K \in \mesh$ (or, equivalently, $\rho_K >0$, $\forall K \in \mesh$).
\end{prpstn}
\begin{proof}
Let us suppose that $\tilde \rho_\edge^{n+1} >0$, $\forall \edge \in \edgesint$.
Then the theory of section \ref{sec:darcy} applies, with:
\[
\normsd{\bar u^{n+1}}=\sum_{\edge \in \edgesint} |D_\edge|\,\frac{\tilde \rho_\edge^{n+1}}{\dt} \, (\bar u_\edge^{n+1})^2
\]
This yields both the existence of a solution and the positivity of the pressure.
In view of the form of the discrete density prediction \eqref{mass-diamond}, this latter property extends by induction to any time step of the computation (provided, of course, that the initial density is positive everywhere).
\end{proof}

\bigskip
We finish this section by some remarks concerning the projection step at hand.

%
\begin{lmm}\label{bmbt}The following identity holds for each discrete pressure $q \in L_h$:
\[
\forall K \in \mesh, \qquad 
(\matB\,\matM_{\tilde \rho^{n+1}}^{-1}\, {\matB}^t\, q)_K = \sum_{\edge=K|L} \frac{1}{\tilde \rho_\edge^{n+1}}\ \frac{|\edge|^2}{|D_\edge|}\, (q_K-q_L)
\]
\end{lmm}

\begin{proof}
Let $q \in L_h$ be given.
By relation \eqref{darcy-disc-deux}, we have:
\[
(\matB^t\,q)_{\edge,i}=|\edge|\ (q_K-q_L)\,n_{KL} \cdot e_i
\]
Let $1_K \in L_h$ be the characteristic function of $K$.
Denoting by $(\cdot,\cdot )$ the standard Euclidian scalar product, by the previous relation and the definition of lumped velocity mass matrix, we obtain:
\[
\begin{array}{ll} \displaystyle
(\matB\,\matM_{\tilde \rho}^{-1}\, \matB^t\, q,\ 1_K)
& \displaystyle
=(\matM_{\tilde \rho}^{-1}\, \matB^t\, q, \ \matB^{t}\,1_K)
\\ & \displaystyle
=\sum_{\edge \in \edgesint} \quad \sum_{i=1}^d \quad \frac{1}{\tilde\rho_{\edge}^{n+1} |D_\edge|}
(\matB^{t}\,q)_{\edge,i} \ (\matB^{t}\,1_K)_{\edge,i}
\\ & \displaystyle
= \sum_{\edge \in \edgesint,\ \edge=K|L} \quad \sum_{i=1}^d \quad \frac{1}{\tilde\rho_{\edge}^{n+1} |D_\edge|}
\left[ |\edge|\, (q_K-q_L)\ n_{KL}\cdot e_i \right] \left[ |\edge|\, n_{KL}\cdot e_i\right]
\end{array}
\]
which, remarking that $\sum_{i=1}^d (n_{KL}\cdot e_i)^2=1$, yields the result.
\end{proof}
\begin{rmrk}[On spurious pressure boundary conditions]
In the context of projection methods for uncompressible flow, it is known that spurious boundary conditions are to be imposed to the pressure in the projection step, in order to make the definition of this step of the algorithm complete.
These boundary conditions are explicit when the process to derive the projection step is first to pose the problem at the time semi-discrete level and then discretize it in space; for instance, with a constant density equal to one and prescribed velocity boundary conditions on $\partial \Omega$, the semi-discrete projection step would take the form:
\[
\left| \begin{array}{ll} \displaystyle
- \Delta\, (p^{n+1}-\tilde p^{n+1})=- \frac 1 \dt \nabla \cdot \tilde u^{n+1}
\hspace{5ex} & \displaystyle
\mbox{in } \Omega
\\[2ex] \displaystyle
\nabla \, (p^{n+1}-\tilde p^{n+1}) \cdot n =0
& \displaystyle
\mbox{on } \partial \Omega
\end{array} \right .
\]
When the elliptic problem for the pressure is built at the algebraic level, the boundary conditions for the pressure are somehow hidden in the discrete operator $\matB\,\matM^{-1}\, {\matB}^t$.
Lemma \ref{bmbt} shows that this matrix takes the form of a finite-volume Laplace discrete operator, with homogeneous Neumann boundary conditions, \ie \ the same boundary conditions as in the time semi-discrete problem above stated.
\label{spurious-BC}\end{rmrk}

\begin{rmrk}[On the non-consistency of the discretization at hand for the Darcy problem]
Considering the semi-discrete problem \eqref{stab-scheme-4}, taking $\tilde \rho_\edge^{n+1}=1,\ \forall \edge \in \edgesint$, one could expect to recover a consistent discretization of a Poisson problem with homogeneous Neumann boundary conditions.
The following example shows that this route is misleading. 
Let us take for the mesh a uniform square grid of step $h$.
The coefficient ${|\edge|^2}/{|D_\edge|}$ can be easily evaluated, and we obtain:
\[
(\matB\,\matM^{-1}\, {\matB}^t\, q)_K = d \sum_{\edge=K|L} \frac{|\edge|}{h}\, (q_K-q_L)
\]
that is the usual finite volume Laplace operator, but multiplied by the space dimension $d$.
This result is of course consistent with the wellknown non-consistency of the Rannacher-Turek element for the Darcy problem, and similar examples could also be given for simplicial grids, with the Crouzeix-Raviart element.
\label{non-consistent!}\end{rmrk}


\subsection{Renormalization steps}\label{subsec:renorm}

The pressure renormalization (step 2 of the algorithm) is introduced for stability reasons, and its form can only be clarified by the stability study.
It reads:
\begin{equation}
\matB\,\matM^{-1}_{\tilde \rho^{n+1}}\, {\matB}^t \ \tilde p^{n+1} = 
\matB\,\matM^{-1}_{\sqrt{\tilde \rho^{n+1}}\,\sqrt{\rho^n}}\, {\matB}^t \ p^n
\label{p-renorm-mat}\end{equation}
where the density at the face at time level $n$ necessary to the definition of $\matM^{-1}_{\sqrt{\tilde \rho^{n+1}}\,\sqrt{\rho^n}}$ is obtained from the density field $\rho^n \in L_h$ by equation \eqref{rho_faces}.
In view of the expression of these operators provided by lemma \ref{bmbt}, this relation equivalently reads:
\begin{equation}
\forall K \in \mesh, \qquad \qquad
\sum_{\edge=K|L} \frac{1}{\tilde \rho_\edge^{n+1}}\ \frac{|\edge|^2}{|D_\edge|}\, (\tilde p_K^{n+1}- \tilde p_L^{n+1})=
\sum_{\edge=K|L} \frac{1}{\sqrt{\tilde \rho_\edge^{n+1}}\sqrt{\rho_\edge^n}}\ \frac{|\edge|^2}{|D_\edge|}\, (p_K^n- p_L^n)
\label{p-renorm}\end{equation}
As $\matB^t$ and $\matB$ stands for respectively the discrete gradient and (opposite of the) divergence operator, this system can be seen as a discretization of the semi-discrete expression of step 2; note however, as detailed in remark \ref{non-consistent!}, that this discretization is non-consistent.

\bigskip
The velocity renormalization (step 5 of the algorithm) simply reads:
\begin{equation}
\forall \edge \in \edgesint, \quad \sqrt{\rho^{n+1}_\edge} \, u^{n+1}_\edge = \sqrt{\tilde \rho^{n+1}_\edge} \, \bar u^{n+1}_\edge
\hspace{10ex} \mbox{or} \qquad
\matM_{\sqrt{\rho^{n+1}}}\ u^{n+1} = \matM_{\sqrt{\tilde \rho^{n+1}}}\ \bar u^{n+1}
\label{v-renorm}\end{equation}


\subsection{An overview of the algorithm}\label{subsec:overview}

To sum up, the algorithm considered in this section is the following one:
\begin{enumerate}
\item Prediction of the density -- The density on the edges at $t^n$, $(\rho_\edge^n)_{\edge \in \edgesint}$, being given by \eqref{rho_faces}, compute $(\tilde \rho^{n+1})_{\edge \in \edgesint}$ by the upwind finite volume discretization of the mass balance over the diamond cells \eqref{mass-diamond}.

\smallskip
\item Renormalization of the pressure -- Compute a renormalized pressure $(\tilde p_K^{n+1})_{K\in\mesh}$ by equation \eqref{p-renorm}.

\smallskip
\item Prediction of the velocity -- Compute $(\tilde u^{n+1}_\edge)_{\edge \in \edgesint}$ by equation \eqref{momentum}, obtained by a finite volume discretization of the transport terms over the diamond cells and a finite element discretization of the other terms.

\smallskip
\item Projection step -- Compute $(\bar u^{n+1}_\edge)_{\edge \in \edgesint}$ and $(p_K^{n+1})_{K\in\mesh}$ from equation \eqref{darcy-disc-un}, obtained by a finite element discretization of the velocity correction equation and an upwind finite volume discretization of the mass balance (over the elements $K \in \mesh$).

\smallskip
\item Renormalization of the velocity -- Compute $(u^{n+1}_\edge)_{\edge \in \edgesint}$ from equation \eqref{v-renorm}.
\end{enumerate}

The existence of a solution to step 4 is proven above; the other problems are linear, and their regularity follows from standard coercivity arguments, using the fact that the discrete densities (\ie\ $\rho^n$ and $\tilde \rho^{n+1}$) are positive, provided that this property is satisfied by the initial condition.


\subsection{Stability analysis} \label{subsec:stab}

In this section, we use the following discrete norm and semi-norm:
\begin{equation}
\begin{array}{ll} \displaystyle
\forall v \in W_h, \qquad
& \displaystyle
\normLdiscd{\tilde \rho}{v} = \sum_{\edge \in \edgesint} |D_\edge|\, \tilde \rho_\edge \, v_\edge^2
\\[2ex]
\forall q \in L_h, \qquad
& \displaystyle
\snormundiscd{\tilde \rho}{q}= \sum_{\edge \in \edgesint,\ \edge=K|L} \frac{1}{\tilde \rho_\edge}\ \frac{|\edge|^2}{|D_\edge|}\, (q_K-q_L)^2
\end{array}
\end{equation}
where $\tilde \rho=(\tilde \rho_\edge)_{\edge \in \edgesint}$ is a family of positive real numbers.
The function $\normLdiscd{\tilde \rho}{\cdot}$ defines a norm over $W_h$, and $\snormundisc{\tilde \rho}{\cdot}$ can be seen as a weighted version of the $H^1$ semi-norm classical in the finite volume context \cite{eym-00-fin}.
The links between this latter semi-norm and the problem at hand are clarified in the following lemma, which is a straightforward consequence of lemma \ref{bmbt}.
%
%
\begin{lmm}\label{bmbt_norm}The following identity holds for each discrete pressure $q \in L_h$:
\[
(\matB\,\matM_{\tilde \rho}^{-1}\, {\matB}^t\, q,q)=\snormundiscd{\tilde \rho}{q}
\]
\end{lmm}

\bigskip
%
%
We are now in position to state the stability of the scheme under consideration.
\begin{thrm}[Stability of the scheme]
Let the equation of state $\varrho(\cdot)$ be such that $\varrho(0)=0$, $\lim_{z \rightarrow + \infty} \varrho(z)= + \infty$ and there exists an elastic potential function $P(\cdot)$ satisfying \eqref{elasticpot_gene} such that the function $z \mapsto z\,P(z)$ is bounded from above in $(0, + \infty)$,  once continuously differentiable and strictly convex.
Let $(\tilde u^n)_{0\leq n \leq N}$, $(u^n)_{0\leq n \leq N}$, $(p^n)_{0\leq n \leq N}$ and $(\rho^n)_{0\leq n \leq N}$ be the solution to the considered scheme, whith a zero forcing term.
Then the following bound holds for $n < N$:
\begin{equation}
\begin{array}{l} \displaystyle
\frac 1 2 \normLdiscd{\rho^{n+1}}{u^{n+1}} + \int_\Omega \rho^{n+1} P(\rho^{n+1}) \dx
+ \dt \sum_{k=1}^{n+1} \int_{\Omega,h} \nabla \tilde u^k:\tau(\tilde u^k)\dx + \frac{\dt^2}{2} \snormundiscd{\tilde \rho^{n+1}}{p^{n+1}}
\\ \displaystyle
\hspace{50ex} \leq \frac 1 2 \normLdiscd{\rho^0}{u^0} + \int_\Omega \rho^0 P(\rho^0) \dx
+ \frac{\dt^2}{2} \snormundiscd{\tilde \rho^0}{p^0}\qquad
\label{stabres}\end{array}
\end{equation}
\label{stab}\end{thrm}

\begin{proof}
Multiplying each equation of the step 3 of the scheme \eqref{momentum} by the corresponding unknown (\textit{i.e} the corresponding component of the velocity $\tilde u^{n+1}$ on the corresponding edge $\sigma$) and summing over the edges and the components yields, by virtue of the stability of the discrete advection operator (theorem \ref{VF1}):
\begin{equation}
\frac{1}{2\, \dt} \normLdiscd{\tilde \rho^{n+1}}{\tilde u^{n+1}} - \frac{1}{2\, \dt} \normLdiscd{\rho^n}{u^n}
+ \int_{\Omega,h} \tau(\tilde u^{n+1}): \nabla \tilde u^{n+1} \dx - \int_{\Omega,h} \tilde p^{n+1} \nabla \cdot \tilde u^{n+1} \dx \leq 0
\label{stab1} \end{equation}
On the other hand, reordering equation \eqref{darcy-alg} and multiplying by $\matM_{\tilde \rho^{n+1}}^{-1/2}$ (recall that $\matM_{\tilde \rho^{n+1}}$ is diagonal), we obtain:
\[
\frac{1}{\dt}\, \matM_{\tilde \rho^{n+1}}^{1/2} \bar u^{n+1} 
+ \matM_{\tilde \rho^{n+1}}^{-1/2}\,\matB^t\, p^{n+1} =
\frac{1}{\dt}\, \matM_{\tilde \rho^{n+1}}^{1/2} \tilde u^{n+1} 
+ \matM_{\tilde \rho^{n+1}}^{-1/2}\,\matB^t\, \tilde p^{n+1}
\]
Squaring this relation gives:
\[
\begin{array}{l}\displaystyle 
\left(\frac{1}{\dt}\, \matM_{\tilde \rho^{n+1}}^{1/2} \bar u^{n+1} 
+ \matM_{\tilde \rho^{n+1}}^{-1/2}\,\matB^t\, p^{n+1},
\ \frac{1}{\dt} \, \matM_{\tilde \rho^{n+1}}^{1/2} \bar u^{n+1} 
+ \matM_{\tilde \rho^{n+1}}^{-1/2}\,\matB^t\, p^{n+1}\right) = 
\\ \displaystyle \hspace{20ex}
\left(\frac{1}{\dt}\, \matM_{\tilde \rho^{n+1}}^{1/2} \tilde u^{n+1} 
+ \matM_{\tilde \rho^{n+1}}^{-1/2}\,\matB^t\, \tilde p^{n+1},
\ \frac{1}{\dt}\, \matM_{\tilde \rho^{n+1}}^{1/2} \tilde u^{n+1} 
+ \matM_{\tilde \rho^{n+1}}^{-1/2}\,\matB^t\, \tilde p^{n+1}\right)
\end{array}
\]
which reads:
\[
\begin{array}{l} \displaystyle 
\frac{1}{\dt^{2}}\, \left( \matM_{\tilde \rho^{n+1}} \bar u^{n+1}, \ \bar u^{n+1}\right) 
+ \left(\matM_{\tilde \rho^{n+1}}^{-1}\,\matB^t\, p^{n+1}, \ \matB^t\, p^{n+1}\right)
+\frac{2}{\dt} \left(\bar u^{n+1},\  \matB^t\, p^{n+1}\right)=
\\ \displaystyle \hspace{20ex}
\frac{1}{\dt^{2}}\, \left(\matM_{\tilde \rho^{n+1}} \tilde u^{n+1}, \ \tilde u^{n+1}\right) 
+\left(\matM_{\tilde \rho^{n+1}}^{-1}\,\matB^t\, \tilde p^{n+1}, \ \matB^t\, \tilde p^{n+1}\right)
+\frac{2}{\dt}\, \left(\tilde u^{n+1}, \ \tilde p^{n+1}\right)
\end{array}
\]
Multiplying by $\dt/2$, remarking that, $\forall v \in W_h,\ (\matM_{\tilde \rho^{n+1}} v, \ v)= \normLdiscd{\tilde \rho^{n+1}}{v}$ and that, thanks to lemma \ref{bmbt_norm}, $\forall q \in L_h,\ (\matM_{\tilde \rho^{n+1}}^{-1}\,\matB^t\, q, \ \matB^t\, q)=
 (\matB \, \matM_{\tilde \rho^{n+1}}^{-1}\,\matB^t\, q,q)=\snormundiscd{\tilde \rho^{n+1}}{q}$, we get: 
\begin{equation}
\begin{array}{l} \displaystyle
\frac{1}{2 \dt} \normLdiscd{\tilde \rho^{n+1}}{\bar u^{n+1}}
+\frac \dt 2 \snormundiscd{\tilde \rho^{n+1}}{p^{n+1}}
+(\bar u^{n+1}, \matB^t\, p^{n+1})
\\ \hspace{20ex} \displaystyle
-\frac{1}{2 \dt} \normLdiscd{\tilde \rho^{n+1}}{\tilde u^{n+1}}
-\frac \dt 2\,\snormundiscd{\tilde \rho^{n+1}}{\tilde p^{n+1}}
- (\tilde u^{n+1}, \matB^t\, \tilde p^{n+1}) =0
\end{array}
\label{stab2}\end{equation}
The quantity $-(\tilde u^{n+1}, \matB^t\, \tilde p^{n+1})$ is nothing more than the opposite of the term  $\displaystyle \int_{\Omega,h} \tilde p^{n+1} \nabla \cdot \tilde u^{n+1} \dx$ appearing in (\ref{stab1}), so summing (\ref{stab1}) and (\ref{stab2}) makes these terms disappear, leading to:
\[
\begin{array}{l} \displaystyle
\frac{1}{2 \dt} \normLdiscd{\tilde \rho^{n+1}}{\bar u^{n+1}}
- \frac{1}{2\, \dt} \normLdiscd{\rho^n}{u^n} 
+ \int_{\Omega,h} \tau(\tilde u^{n+1}) : \nabla \tilde u^{n+1} \dx
\\[2ex] \hspace{20ex} \displaystyle
+\frac \dt 2 \, \snormundiscd{\tilde \rho^{n+1}}{p^{n+1}}
-\frac \dt 2 \, \snormundiscd{\tilde \rho^{n+1}}{\tilde p^{n+1}}
+ (\bar u^{n+1}, \matB^t\, p^{n+1})\leq 0 
\end{array}
\]
Finally, $(\bar u^{n+1}, \matB^t\, p^{n+1})$ is precisely the pressure work which can be bounded by the time derivative of the elastic potential, as stated in theorem \ref{VF2}:
\begin{equation}
\begin{array}{l} \displaystyle
\frac{1}{2 \dt} \normLdiscd{\tilde \rho^{n+1}}{\bar u^{n+1}} 
+ \int_{\Omega,h} \tau(\tilde u^{n+1}) : \nabla \tilde u^{n+1} \dx 
+\frac \dt 2 \,\snormundiscd{\tilde \rho^{n+1}}{p^{n+1}}
-\frac \dt 2 \,\snormundiscd{\tilde \rho^{n+1}}{\tilde p^{n+1}} 
\\ \hspace{30ex} \displaystyle
+\frac{1}{\dt}  \int_\Omega \rho^{n+1} P(\rho^{n+1}) \dx
\leq \frac{1}{2\, \dt} \normLdiscd{\rho^n}{u^n} 
+\frac{1}{\dt} \int_\Omega \rho^{n} P(\rho^{n}) \dx
\end{array}\label{eqstep2}\end{equation}
The proof then ends by using the renormalization steps (step 2 and 5 of the algorithm). 
Step 2 reads in an algebric setting:
\[
\matB\, \matM_{\tilde \rho^{n+1}}^{-1}\,\matB^t\, \tilde p^{n+1}
=\matB\, \matM_{\tilde \rho^{n+1}}^{-1/2}\, \matM_{\tilde \rho^{n}}^{-1/2}\,\matB^t\,  p^{n}
\]
Multiplying by $\tilde p^{n+1}$, we obtain:
\[
\left(\matM_{\tilde \rho^{n+1}}^{-1/2}\,\matB^t\, \tilde p^{n+1},\ \matM_{\tilde \rho^{n+1}}^{-1/2}\,\matB^t\, \tilde p^{n+1}\right)
= \left(\matM_{\tilde \rho^{n}}^{-1/2}\,\matB^t\, p^{n}, \ \matM_{\tilde \rho^{n+1}}^{-1/2}\, \matB^t\, \tilde p^{n+1}\right)
\]
and thus, by Cauchy-Schwartz inequality:
\[
\left(\matM_{\tilde \rho^{n+1}}^{-1/2}\,\matB^t\, \tilde p^{n+1},\ \matM_{\tilde \rho^{n+1}}^{-1/2}\,\matB^t\, \tilde p^{n+1}\right)
\leq
\left(\matM_{\tilde \rho^{n}}^{-1/2}\,\matB^t\, p^{n}, \ \matM_{\tilde \rho^{n}}^{-1/2}\,\matB^t\, p^{n}\right)^{1/2}\, 
\left(\matM_{\tilde \rho^{n+1}}^{-1/2}\,\matB^t\, \tilde p^{n+1},\ \matM_{\tilde \rho^{n+1}}^{-1/2}\,\matB^t\, \tilde p^{n+1} \right)^{1/2}
\]
This relation yields $\snormundiscd{\tilde \rho^{n+1}}{\tilde p^{n+1}} \leq \snormundiscd{\tilde \rho^{n}}{p^{n}}$.
In addition, step 5 of the algorithm gives $\normLdiscd{\rho^{n+1}}{u^{n+1}}=\normLdiscd{\tilde \rho^{n+1}}{\bar u^{n+1}}$.
Using these two relations in \eqref{eqstep2}, we get:
\[
\begin{array}{l} \displaystyle
\frac{1}{2 \dt} \normLdiscd{\rho^{n+1}}{u^{n+1}} 
+ \int_{\Omega,h} \tau(\tilde u^{n+1}): \nabla \tilde u^{n+1} \dx 
+\frac \dt 2 \,\snormundiscd{\tilde \rho^{n+1}}{p^{n+1}} 
+ \frac{1}{ \dt} \int_\Omega \rho^{n+1} P(\rho^{n+1}) \dx 
\\[2ex] \hspace{50ex} \displaystyle
\leq \frac{1}{2\, \dt} \normLdiscd{\rho^n}{u^n} 
+ \frac{1}{\dt} \int_\Omega \rho^{n} P(\rho^{n}) \dx
+ \frac \dt 2\,\snormundiscd{\tilde \rho^{n}}{p^n}
\end{array}
\]
and the estimate of theorem \ref{stab} follows by summing over the time steps.
\end{proof}
%
%
\begin{rmrk}[On the upwinding of the mass balance discretization, the \textit{inf-sup} stability of the discretization and the appearance of spurious pressure wiggles.]
In the scheme considered in this section, the upwinding in the discretization of mass balance controls the onset of density oscillations.
As long as the pressure and the density are linked by an increasing function, that is as long as the flow remains compressible with a reasonable equation of state, it will probably be sufficient to prevent from pressure oscillations.
Besides, the fourth term of the left hand side of (\ref{stabres}), \textit{i.e.} the term involving $\snormundiscd{p^{n+1}}{\tilde \rho^{n+1}}$,  provides a control on the discrete $H^1$ semi-norm of the pressure, at least for large time steps, and therefore also produces an additional pressure smearing.
However, it comes up in the analysis as the composition of the discrete divergence with the discrete gradient; consequently, one will obtain such a smoothing effect only for \textit{inf-sup} stable discretizations.
Note also that, even for steady state problems, some authors recommend the use of stable approximation space pairs to avoid pressure wiggles \cite{bris-90-com, for-93-fin}.
\end{rmrk}

\bigskip
\begin{rmrk}[On a different projection step]
Some authors propose a different projection step \cite{bij-98-uni,wes-01-pri}, which reads in time semi-discrete setting:
\[
\left| \begin{array}{l} \displaystyle
\frac{\varrho(p^{n+1})\, u^{n+1}- \tilde \rho^{n+1}\, \tilde u^{n+1}}{\dt}
+ \nabla (p^{n+1}-\tilde p^{n+1})=0
\\[2ex] \displaystyle
\frac{\varrho(p^{n+1})-\rho^n}{\dt}
+ \nabla \cdot \left( \varrho(p^{n+1})\, u^{n+1} \right)=0
\end{array} \right .
\]
Considering this system, one may be tempted by the following line of thought: choosing $q^{n+1}=\varrho(p^{n+1})\, u^{n+1}$ as variable, taking the discrete divergence of the first equation and using the second one will cause the convection term of the mass balance to  disappear from the discrete elliptic problem for the pressure, whatever the discretization of this term may be.
Consequently, the equation for the pressure will be free of the non-linearities induced by the upwinding and the dependency of the convected density on the pressure, while one still may hope to obtain a positive upwind (with respect to the density) scheme.
In fact, this last point is incorrect.
To be valid, it would necessitate that, from any solution $(q^{n+1}, p^{n+1})$, one be able to compute a velocity field $u^{n+1}$ by dividing $q^{n+1}$ by the density of the control volume located upstream with respect to $u^{n+1}$.
Unfortunately, it is not always possible to obtain this upstream value; for instance, if for two neighbouring control volumes $K$ and $L$, $\rho_K<0$, $\rho_L>0$ and $q^{n+1} \cdot n_{K|L}>0$, neither the choice of $K$ nor $L$ for the upstream control volume is valid.
Consequently, with this discretization, we are no longer able to guarantee neither the positivity of $\rho$ nor the absence of oscillations.
However, as explained below, if the density remains positive, we will have a smearing of pressure or density wiggles due to the fact that the discretization is \textit{inf-sup} stable.
\end{rmrk}


\subsection{Implementation}\label{subsec:impl}

The implementation of the first three steps (\ref{stab-scheme-1})-(\ref{stab-scheme-3}) of the numerical scheme  is standard, and we therefore only describe here in details the fourth step, that is the projection step.
The precise algebraic formulation of the system (\ref{stab-scheme-4}) reads:
\begin{equation}
\begin{array}{l} 
\left| \begin{array}{l} \displaystyle
\frac{1}{\dt} \matM_{\tilde \rho^{n+1}} \, (\bar u^{n+1} - \tilde u^{n+1})
+ \matB^t\, (p^{n+1} -\tilde p^{n+1})=0
\\[2ex] \displaystyle
\frac{1}{\dt} \matR\, (\varrho(p^{n+1})-\rho^n) - \matB \matQ_{\rho^{n+1}_{\rm up}} \bar u^{n+1}=0
\end{array} \right .
\end{array}
\label{proj-alg}\end{equation}
where ${\rm M}_{\tilde \rho^{n+1}}$ and ${\rm Q}_{\rho^{n+1}_{\rm up}}$ are two diagonal matrices; for the first one, we recall that the entry corresponding to an edge $\edge\in \edgesint,\ \edge=K|L$ is computed by multiplying the measure of the diamond associated to $\edge$ by the predicted density (at the edge center) $\tilde \rho^{n+1}_\edge$; in the second one, the same entry is obtained by just taking the density at $t^{n+1}$ in the element located upstream of $\edge$ with respect to $\bar u^{n+1}$, {\it i.e.} either $\varrho\,(p_K^{n+1})$ or $\varrho\,(p_L^{n+1})$.
Note that these definitions can be extended in a straightforward way for the boundary edges, if the velocity is not prescribed to zero on the boundary of the computational domain.
The matrix ${\matR}$ is diagonal and, for any $K\in\mesh$, its entry $\matR_K$ is the measure of the element $K$.
For the sake of simplicity, we suppose for the moment that the equation of state is linear:
\[
\varrho(p^{n+1})=\frac{\partial \varrho}{\partial p}\ p^{n+1}
\]

\medskip
The elliptic problem for the pressure is obtained by multiplying the first relation of (\ref{proj-alg}) by $\matB\ \matQ_{\rho^{n+1}_{\rm up}}\ (\matM_{\tilde \rho^{n+1}})^{-1}$ and using the second one.
This equation reads:
\begin{equation}
\matL\, p^{n+1}+\frac{\partial \varrho}{\partial p}\ \frac{1}{\dt^{2}}\matR\, p^{n+1} =
\matL\,\tilde p^{n+1}+\frac{1}{\dt^{2}}\matR \rho^n +\frac{1}{\dt}\matB\,\matQ_{\rho^{n+1}_{\rm up}}\,\tilde u^{n+1}
\label{step41}\end{equation}
where $\matL=\matB\,\matQ_{\rho^{n+1}_{\rm up}}\,(\matM_{\tilde \rho^{n+1}})^{-1}\,\matB^t$ can be viewed, for the discretization at hand, as a finite volume discrete approximation of the Laplace operator with Neumann boundary conditions (when the velocity is prescribed at the boundary), multiplied by the space dimension $d$ and the densities ratio (see remarks \ref{spurious-BC} and \ref{non-consistent!}).
We recall that, by a calculation similar to the proof of lemma \ref{bmbt}, this matrix can be evaluated directly in the "finite volume way", by the following relation, valid for each element $K$:
\[
(\matL\, p^{n+1})_K=\sum_{\edge=K|L} \ 
d\ \frac{\displaystyle \rho_{{\rm up},\edge}}{\tilde \rho^{n+1}_\edge}\ \frac{|\edge|^2}{|D_\edge|}\, (p_K-p_L)
\]
where $\rho_{{\rm up},\edge}$ stands for the upwind density associated to the edge $\edge$.
Provided that $p^{n+1}$ is known, the first relation of (\ref{proj-alg}) gives us the updated value of the velocity:
\[
\bar u^{n+1}=\tilde u^{n+1}-\dt\ (\matM_{\tilde \rho^{n+1}})^{-1}\,\matB^t\,(p^{n+1}-\tilde p^{n+1})
\label{step42}
\]
As, to preserve the positivity of the density, we want to use in the mass balance the value of the density upwinded with respect to $\bar u^{n+1}$, equations (\ref{step41}) and (\ref{step42}) are not decoupled, by contrast with what happens in usual projection methods.
We thus implement the following iterative algorithm:
\[
\begin{array}{ll} \displaystyle
\mbox{Initialization: }\quad
p_{0}^{n+1}= \tilde p^{n+1} \mbox{ and } \bar u_{0}^{n+1}=\tilde u^{n+1}
\\[3ex] 
\qquad \left| \begin{array}{l} \displaystyle
\mbox{Step 4.1 -- Solve for } p_{k+1/2}^{n+1}\ :
\\[3ex] \displaystyle \hspace{5ex}
\left[\matL+\frac{\partial \varrho}{\partial p}\ \frac{1}{\dt^{2}} \matR \right]\, p_{k+1/2}^{n+1}=
\matL\,\tilde p^{n+1} + \frac{1}{\dt^{2}} \matR \,  \rho^n
+\frac{1}{\dt} \matB\, \matQ_{\rho^{n+1}_{\rm up}}\, \tilde u^{n+1}
\\[3ex] \displaystyle
\quad \mbox{where the density in $\matL$ and $\matQ_{\rho^{n+1}_{\rm up}}$ is evaluated at $p^{n+1}_k$ and the upwinding}
\\ \displaystyle 
\quad \mbox{in $\matQ_{\rho^{n+1}_{\rm up}}$ is performed with respect to $\bar u^{n+1}_k$}
\\[3ex] \displaystyle
\mbox{Step 4.2 -- Compute } p_{k+1}^{n+1}\mbox{ as } p_{k+1}^{n+1}=\alpha \, p_{k+1/2}^{n+1} + (1-\alpha)\, p_{k}^{n+1}
\\[3ex] \displaystyle
\mbox{Step 4.3 -- Compute } \bar u_{k+1}^{n+1}\mbox{ as }:
\\[2ex] \displaystyle \hspace{5ex}
\bar u_{k+1}^{n+1}=\tilde u^{n+1}-\dt \ (\matM_{\tilde \rho^{n+1}})^{-1}\,\matB^t\,(p_{k+1}^{n+1}-\tilde p^{n+1})
\end{array} \right.
\\[20ex] \displaystyle 
\mbox{Convergence criteria: }\quad 
\max \left[ \normz{p_{k+1}^{n+1}-p_{k}^{n+1}},\normz{u^{n+1}_{k+1}-u_{k}^{n+1}}\right] < \varepsilon
\end{array}
\]
The second step of the previous algorithm is a relaxation step which can be performed to ensure convergence; however, in the tests presented hereafter, we use $\alpha=1$ and obtain convergence in few iterations (typically less than 5).
When the equation of state is nonlinear, step 4.1 is replaced by one iteration of Newton's algorithm.


\subsection{Numerical experiments}\label{subsec:num-exp}
In this section, we describe numerical experiments performed to assess the behaviour of the pressure correction scheme presented in this paper, in particular the convergence rate with respect to the space and time discretizations.

\medskip
With $\Omega=(0,1)\times (-\frac 1 2,\frac 1 2)$, we choose for the momentum and density the following expressions:
\[
\begin{array}{l} \displaystyle
\rho\,u= -\frac{1}{4} \cos(\pi t)\left[\begin{array}{l} \sin(\pi x^{(1)}) \\ \cos(\pi x^{(2)}) \end{array}\right]
\\[3ex] \displaystyle
\rho=1+\frac{1}{4}\,\sin(\pi t)\,\left[\cos(\pi x^{(1)})-\sin(\pi x^{(2)})\right]
\end{array}
\]
These functions satisfy the mass balance equation; for the momentum balance, we add the corresponding right hand side.
In this latter equation, the divergence of the stress tensor is given by:
\[
\nabla \cdot \tau(u)= \mu \Delta u + \frac \mu 3 \nabla\, \nabla \cdot u,
\qquad \mu=10^{-2}
\]
The pressure is linked to the density by the following equation of state:
\[
p=\varrho(\rho)=\frac{\rho - 1}{\gamma\ {\rm Ma}^2}, \qquad \gamma = 1.4, \ {\rm Ma}=0.5
\]
where the parameter $\rm Ma$ corresponds to the characteristic Mach number.

\begin{figure}
\begin{center} \scalebox{0.6}{\includegraphics*{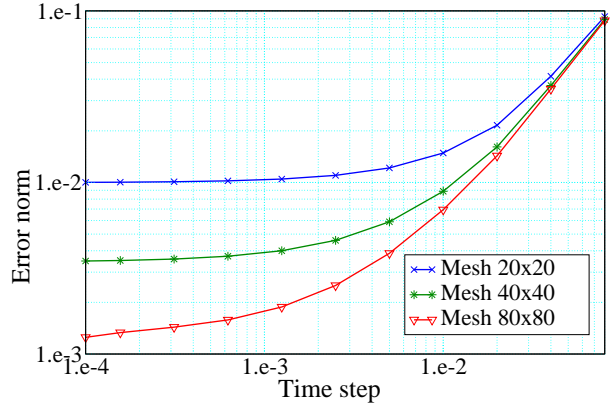}} \end{center}
\caption{velocity error as a function of the time step.
\label{err_v}}
\end{figure}

\medskip
We use in these tests a special numerical integration of the forcing term of the momentum balance, which is designed to ensure that the discretization of a gradient is indeed a discrete gradient ({\it i.e.} if the forcing term $f$ can be recast under the form $f=\nabla g$, the discrete right hand side of the momentum balance belongs to the range of ${\rm B^t}$).

\begin{figure}
\begin{center} \scalebox{0.6}{\includegraphics*{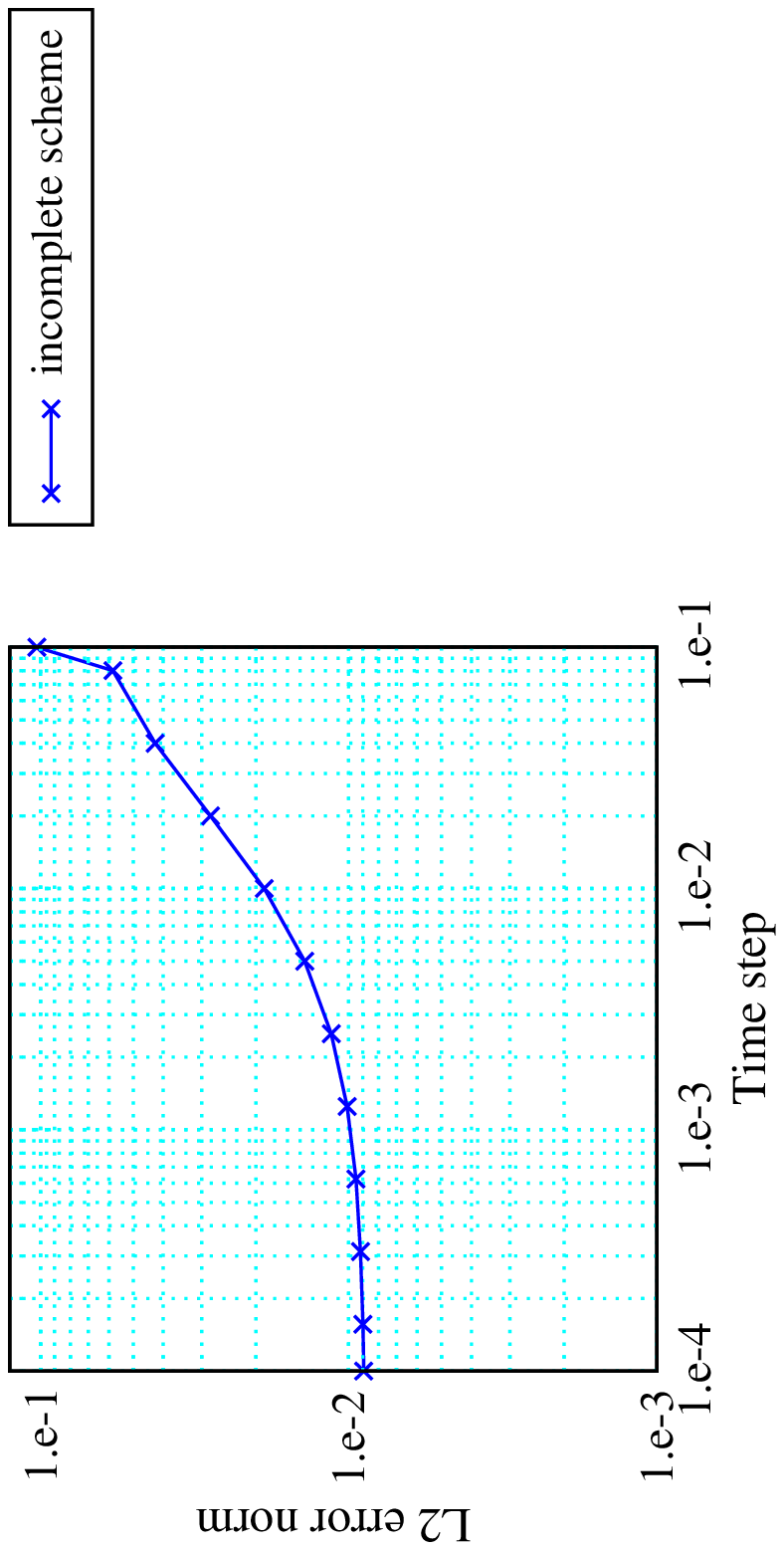}} \end{center}
\caption{pressure error as a function of the time step.
\label{err_p}}
\end{figure}

\medskip
Velocity and pressure errors obtained at $t=0.5$, in respectively $L^2$ and discrete $L^2$ norms and as a function of the time step, are drawn on respectively figure \ref{err_v} and figure \ref{err_p}, for $20 \times 20$, $40 \times 40$ and $80 \times 80$ uniform meshes.
For large time steps, these curves show a decrease corresponding to approximately a first order convergence in time for the velocity and the pressure, until a plateau is reached, due to the fact that errors are bounded from below by the residual spatial discretization error.
For both velocity and pressure, the value of the errors on this plateau show a spatial convergence order (in $L^2$ norm) close to 2.


\section{Conclusion}

We presented in this paper a numerical scheme for the barotropic Navier-Stokes compressible equations, based on a pressure-correction time stepping algorithm.
For the spatial discretization, it combines low-order non-conforming mixed finite elements with finite volumes; in the incompressible limit, one recovers a classical projection scheme based on an \textit{inf-sup} stable pair of approximation spaces for the velocity and the pressure.
This scheme is proven to enjoy an unconditional stability property: irrespectively of the time step, the discrete solution obeys the \textit{a priori} estimates associated to the continuous problem, \textit{i.e.} strict positivity of the density, bounds in $L^\infty$-in-time norm of the quantity $\int_\Omega \rho\, u^2\, {\rm d}x$ and $\int_\Omega \rho \, P(\rho)\, {\rm d}x$ and in $L^2$-in-time norm of the viscous dissipation $ \int_\Omega \tau(u): \nabla u \, {\rm d}x$.
To our knowledge, this result is the first one of this type for barotropic compressible flows.

\medskip
However, the scheme presented here is by no means "the ultimate scheme" for the solution to the compressible Navier-Stokes equations.
It should rather be seen as an example of application (and probably one of the less sophisticated ones) of the mathematical arguments developped to obtain  stability, namely theorems \ref{VF2} (discrete elastic potential identity) and \ref{VF1} (stability of the advection operator), and our hope is that these two ingredients could be used as such or adapted in the future to study other algorithms.  
For instance, a computation close to the proof of theorem \ref{stab} (and even simpler) would yield the stability of the fully implicit scheme; adding to this latter algorithm a prediction step for the density (as performed here) would also allow to linearize (once again as performed here) the convection operator without loss of stability.
A stable pressure-correction scheme avoiding this prediction step can also be obtained, and is currently under tests at IRSN for the computation of compressible bubbly flows.
Besides these variants, less diffusive schemes should certainly be searched for.
Finally, the proposed scheme is currently the object of more in-depth numerical studies including, in particular, problems admitting less smooth solutions than the test presented here.


\bibliographystyle{plain}
\bibliography{./nscomp}
\end{document}